\documentclass[12pt]{article}
\usepackage{amssymb}
\usepackage{a4}
\begin{document}
\def\st{\, : \,}
\def\kbar{{\mathchar'26\mkern-9muk}}  
\def\bra#1{\langle #1 \vert}
\def\ket#1{\vert #1 \rangle}
\def\vev#1{\langle #1 \rangle}
\def\ad{\mbox{ad}\,}
\def\ker{\mbox{Ker}\,}
\def\im{\mbox{Im}\,}
\def\der{\mbox{Der}\,}
\def\ad{\mbox{ad}\,}
\def\b#1{{\mathbb #1}}
\def\c#1{{\cal #1}}
\def\pt{\partial_t}
\def\px{\partial_1}
\def\bpx{\bar\partial_1}
\def\la{\langle}
\def\ra{\rangle}
\def\nn{\nonumber \\}
\def\pn{\par\noindent}
\def\etal{{\it et al.}\ }
\def\1{{\bf 1}}
\newcommand{\botimes}{\mbox{\boldmath \boldmath$\otimes$}}

\def\R{{\cal R}\,}
\newcommand{\tr}{\triangleright\,}
\newcommand{\tl}{\,\triangleleft}
\newcommand{\tro}{\triangleright^{op}\,}
\newcommand{\tlo}{\,\stackrel{op}{\triangleleft}}
\def\cross{{\triangleright\!\!\!<}}
\def\cocross{{>\!\!\!\triangleleft\,}}
\def\g{\mbox{\bf g\,}}
\def\id{\mbox{id\,}}
\def \uqg{\mbox{$U_q{\/\mbox{\bf g}}$ }}
\def \uqgo{\mbox{$U^{op}_q{\/\mbox{\bf g}}$ }}
\def \uqgp{\mbox{$U_q^+{\/\mbox{\bf g}}$ }}
\def \uqgn{\mbox{$U_q^-{\/\mbox{\bf g}}$ }}

\renewcommand{\thefootnote}{\fnsymbol{footnote}}

\newcommand{\sect}[1]{\setcounter{equation}{0}\section{#1}}
\renewcommand{\theequation}{\thesection.\arabic{equation}}
\newcommand{\subsect}[1]{\setcounter{equation}{0}\subsection{#1}}
\renewcommand{\theequation}{\thesection.\arabic{equation}}
\newcommand{\app}[1]{\setcounter{section}{0}
\setcounter{equation}{0} \renewcommand{\thesection}
{\Alph{section}}\section{#1}}

\newcommand{\be}{\begin{equation}}
\newcommand{\ee}{\end{equation}}
\newcommand{\ba}{\begin{eqnarray}}
\newcommand{\ea}{\end{eqnarray}}
%
%
%
\newtheorem{prop}{Proposition}
\newtheorem{lemma}{Lemma}
\newtheorem{theorem}{Theorem}
\newtheorem{corollary}{Corollary}
%
%
%
%
\def\sq{\mbox{\rlap{$\sqcap$}$\sqcup$}}
\newenvironment{proof}[1]{\vspace{5pt}\noindent{\bf Proof #1}\hspace{6pt}}%
{\hfill\sq}
\newcommand{\bp}{\begin{proof}}
\newcommand{\ep}{\end{proof}\par\vspace{10pt}\noindent}
%
%

\title{On the Decoupling of the Homogeneous and Inhomogeneous Parts
in Inhomogeneous Quantum Groups}

\author{        Gaetano Fiore, \\\\
         \and
        Dip. di Matematica e Applicazioni, Fac.  di Ingegneria\\ 
        Universit\`a di Napoli, V. Claudio 21, 80125 Napoli
        \and
        I.N.F.N., Sezione di Napoli,\\
        Complesso MSA, V. Cintia, 80126 Napoli
        }
\date{}

\maketitle
\abstract{We show that, if there exists a realization of a Hopf
algebra $H$ in a $H$-module algebra $\c{A}$, then one can split
their cross-product into the tensor product algebra of $\c{A}$ itself
with a subalgebra isomorphic to $H$ and {\it commuting} with 
$\c{A}$. This result applies in particular to the
algebra underlying inhomogeneous quantum groups like 
the Euclidean ones, which are obtained as cross-products of the quantum
Euclidean spaces $\b{R}_q^N$ with the quantum groups of rotation 
$U_qso(N)$ of $\b{R}_q^N$, for which it has no classical analog.
}

\vfill
\noindent
Preprint 00-31 Dip. Matematica e Applicazioni, Universit\`a di Napoli\\
DSF/3-2001
\newpage

\sect{Introduction}

As known, given a unital module algebra $\c{A}$ of a Lie algebra
$\g$ (over the field $\b{C}$, say), one can build a new module algebra, 
called cross-product $U\g \cross \c{A}$,
that is as a vector space the tensor product 
$\c{A}\otimes U\g$ of the vector spaces $\c{A}, U\g$ (over the same
field) and has the product law
\ba
&&(\1_{\c{A}}\otimes g)(\1_{\c{A}}\otimes g')= (\1_{\c{A}}\otimes gg'), 
\label{trivial1} \\
&&(a\otimes \1_H)(a'\otimes\1_H)= 
aa'\otimes \1_H, \label{trivial2}\\
&&(a\otimes \1_H)(\1_{\c{A}}\otimes g)= a\otimes g, \label{trivial3} \\
&&(\1_{\c{A}}\otimes g)(a\otimes\1_H)=g_{(1)}\tr a \otimes g_{(2)},
\label{trivial4}
\ea
for any $g,g'\in U\g$, $a,a'\in\c{A}$.
Here we have denoted by $\tr$ the left action of
the Hopf algebra $H\equiv U\g$ on $\c{A}$,
\be
\tr: \:(g,a)\in H\times \c{A} \rightarrow g\tr \!a\in\c{A}, \label{act1}
\ee
and used a Sweedler-type notation
with suppressed summation sign for the coproduct $\Delta(g)$ 
of $g$, namely the short-hand notation 
$\Delta(g)=g_{(1)}\botimes g_{(2)}$ instead of a sum 
$\Delta(g)=\sum_{\mu}g_{(1)}^{\mu}\botimes g_{(2)}^{\mu}$
of many terms. 
In the main part of this paper we shall work with a left action and
therefore left-module algebras. In section \ref{rightsec}
we shall give the formulae if we use instead a right action $\tl$ and
right-module algebras. By definition of a (left) action,
for any $g,g'\in H$, $a,a'\in\c{A}$,
\ba
&& (gg')\tr a=     g\tr(g'\tr a)              \label{modalg1r}\\
&& g\tr(aa')=(g_{(1)}\tr a)\, (g_{(2)}\tr a');       \label{modalg2r}
\ea
We recall that
the coproduct on the unit and on any $g\in\g$ is given by
\be
\Delta(\1_H)=\1_H\otimes\1_H \qquad\qquad
\Delta(g)=g\otimes\1_H+\1_H \otimes g;
\ee
on the rest of $U\g$ it is determined using the fact
that it is an algebra homomorphism
$\Delta: H\to H\botimes H$. Clearly, the
coproduct is cocommutative, i.e.
$g_{(1)}\botimes g_{(2)}=g_{(2)}\botimes g_{(1)}$. 

For the sake of clearness, above and in the sequel we denote
the tensor products
of vector spaces and of algebras by $\otimes$ and $\botimes$
(in boldface) respectively; therefore, given two unital
algebras $B,B'$ (over the same field), $B\botimes B'$ 
is $B\otimes B'$ as a vector space, while as an algebra it
is characterized by the product
\be
(a \otimes a')(b\otimes b')= (ab\otimes a'b')      \label{trivial4'}
\ee
for any $a,b\in B$ and $a',b'\in B'$.

In the sequel, with a standard abuse of notation, 
for any unital algebras $B,B'$ and any $b\in B$, $b'\in B'$
we shall denote by $bb'$ the element $b \otimes b'$ in the
tensor product of vector spaces $B\otimes B'$ and omit either unit
$\1_B,\1_{B'}$ 
whenever multiplied by non-unit elements (thus,
$\1_BB',B\1_{B'}$ will be denoted by $B',B$).
Consequently, in the case
of e.g. the cross product algebra $U\g \cross \c{A}$ 
relations (\ref{trivial1}-\ref{trivial3})
take trivial forms, whereas (\ref{trivial4}) becomes the
commutation relation
\be
g a = (g_{(1)}\tr a) g_{(2)};                    \label{gxrel}
\ee
whereas for a tensor product algebra $B\botimes B'$ the analogs
of relations (\ref{trivial1}-\ref{trivial3}) again
take trivial forms, whereas (\ref{trivial4'}) for 
$a=\1_B$, $b'=\1_{B'}$ becomes
the trivial commutation relation
\be
a' \,b= b\,a'.                                    \label{trivial}
\ee
Of course, $B\botimes B'$
is isomorphic to $B'\botimes B$.

$H \cross \c{A}$ is itself a module algebra under the left
action  $\tr$ of $ H$ if we extend the latter on the elements of 
the $H$ subalgebra as the adjoint action,
\be
g\tr h=g_{(1)} h Sg_{(2)}    \qquad\qquad g,h\in  H; \label{adjo}
\ee
(here $S$ denotes the antipode of $ H$), and set as usual
$$
g\tr(a\otimes h)=g_{(1)}\tr a\:\otimes\: g_{(2)}\tr h;
$$
note that, in the notation mentioned above this relation takes
the same form as (\ref{modalg2r}), i.e. becomes
$g\tr(ah)=(g_{(1)}\tr a)(g_{(2)}\tr h)$.
It is immediate to show that relation (\ref{gxrel}) implies that
one can realize the action 
$\tr:H\times(H \cross\c{A})\rightarrow H \cross\c{A}$ 
in the `adjoint-like way'
\be
g\tr\eta=g_{(1)} \, \eta \, Sg_{(2)}               \label{realiz1}
\ee
on all of $H \cross\c{A}$. 

Classical examples of cross product algebras
are the universal enveloping algebras of inhomogeneous
Lie groups, like the Poincar\'e algebra, where
$H=Uso(3,1)$ is the algebra generated by infinitesimal Lorentz
transformations and $\c{A}$ is the abelian algebra generated by
infinitesimal translations on Minkowski space, or the Euclidean
algebra, where $\c{A}$ is the abelian algebra generated by
infinitesimal translations on the $N$-dimension\-al Euclidean space 
$\b{R}^N$ and
$H=Uso(N)$ is the algebra generated by its infinitesimal rotations.

The above setting can be ``deformed'' by allowing $H$ to be
a non-cocommutative Hopf algebra, e.g. the quantum group $\uqg$, and as
$\c{A}$ the corresponding $q$-deformed module algebra. In general, 
the latter will be no more abelian, even if its classical 
counterpart is. A cross product
$H \cross \c{A}$ can be still defined by means of the same formulae
(\ref{act1}-\ref{gxrel}).

\medskip
In this paper we want to show that 
there are prominent examples of cross products
$H \cross \c{A}$ that are isomorphic to  $\c{A}\botimes H$, 
more precisely are equal to $\c{A}\, H'$
with $H'\subset H \cross \c{A}$ 
a subalgebra isomorphic to $H$ and {\it commuting} 
with $\c{A}$, even if this is not the case
for their undeformed counterparts. 
As we shall see, this occurs if there 
exists an algebra homomorphism $\varphi$
of the cross product into $\c{A}$ acting identically
on the latter.
Of course, this will have dramatic consequences for the cross product
both from the algebraic and from the representation-theoretic viewpoint;
it will allow to reduce representations of $H \cross \c{A}$ to direct sums
of tensor products of representations of $\c{A}$ and of $H'$.
To prevent misunderstandings, we note that in general if $H \cross \c{A}$
itself is a Hopf algebra, as in the case of inhomogeneous quantum 
groups, neither $\c{A}$ nor $H'$ will be a Hopf subalgebra of
$H \cross \c{A}$.

The present work has been inspired among others by the results of Ref.
\cite{CerMadSchWes00,Fio95}. In Ref. \cite{CerMadSchWes00}
the existence of such a $H'$ for the $q$-deformed Euclidean algebra in three
dimensions has been noted; its generators
have been constructed ``by hand'' and have been used to decouple
$q$-rotations from $q$-translations in the
$*$-representations of $H \cross \c{A}$. In  Ref. \cite{Fio95}
we had constructed ``by hand''
for the $q$-deformed Euclidean algebra in $N\ge 3$
dimensions a set of generators which do the same job; but, instead
of commuting with the $q$-deformed generators of translations,
they $q$-commute with the latter. 
Now, {\it a posteriori}, one can check
that they can be obtained as products of suitable
elements of $H'$ by
suitable elements of the natural Cartan (i.e. maximal
abelian) Hopf subalgebra of $H=U_qso(N)$.

On the contrary, the present work gives 
a very simple prescription for their
construction, based on the existence of $\varphi$. 
The prescription can be thus applied to
a number of models, including the following.
In Ref. \cite{CerFioMad00,CerMadSchWes00} a class of
homomorphisms (\ref{Hom}-\ref{ident0})
has been determined for a slightly enlarged version $\c{A}$ of 
the algebra of functions on the $N$-dimensional quantum 
Euclidean space $\b{R}_q^N$ or quantum Euclidean sphere $S_q^{N\!-\!1}$, 
the Hopf algebra $H$ denoting 
$U_qso(N)$ itself if $N$ is odd, either the Borel subalgebra
$U^+_qso(N)$ or the one $U^-_qso(N)$ if $N$ is even.
Their behaviour under the
$*$-structures has been investigated in Ref. \cite{FioSteWes00}
(where incidentally we draw another consequence of its existence, namely
the possibility of ``unbraiding'' braided tensor product algebras).
This will be explicitly described in section \ref{appliEuc}.
The analogous maps for the $q$-deformed fuzzy sphere $S_{q,M}^2$ 
have been found in \cite{fuzzyq}. On the other hand, the existence 
of algebra homomorphisms  (\ref{Hom}) for 
$H=U_qso(N),U_qsl(N)$  and $\c{A}$ respectively equal to
(a suitable completion of) the
$U_qso(N)$-covariant Heisenberg algebra or the $U_qsl(N)$-covariant 
Heisenberg (or Clifford) algebras,
has been known for even a longer time \cite{Fiocmp95,ChuZum95,Hay90}.
This will be treated in section \ref{heisenberg}.
Note that in the latter cases $\varphi$ exist also for the
undeformed counterparts at $q=1$, thus our results will apply 
also in this case (we do not know whether this has ever been 
formulated as a result in ordinary Lie group theory).

In section \ref{commutant} we state and prove
the main results of this work leading to the construction of $H'$.
In section \ref{star} we focus on the $*$-structures.
In section \ref{rightsec} we give without proof for right-module 
algebras all the main formulae valid for left-module algebras.
In section \ref{applications} we apply our results to the two 
examples of cross-product algebras mentioned above.
\medskip

We conclude this section with some additional preliminaries.
Beside the Sweedler-type notation with lower indices introduced
for the coproduct, we shall denote a sum of many terms in a tensor 
product by a Sweedler-type notation with {\it upper} indices
and suppressed summation sign, 
e.g. $c^{(1)}\otimes c^{(2)}$ will actually mean a sum
$\sum_{\mu}c^{(1)}_{\mu}\otimes c^{(2)}_{\mu}$. 
Secondly, we can also introduce an `opposite' action 
$\tro:H\times(H \cross\c{A})\rightarrow H \cross\c{A}$,
i.e. an action of the Hopf algebra with the same algebra structure
and counit but opposite coproduct $\Delta^{op}(g)=g_{(2)}\botimes g_{(1)}$
and inverse antipode,
by
\be
g\tro\eta=g_{(2)} \, \eta \, S^{-1}g_{(1)}.                \label{realizop}
\ee
It fulfills
\ba
&& (gg')\tro a=     g\tro(g'\tro a)              \label{modalg1op}\\
&& g\tro(aa')=(g_{(2)}\tro a)\, (g_{(1)}\tro a').       \label{modalg2op}
\ea
(Note that $g\tro a$ in general is {\it not} an element of $\c{A}$).

\sect{The commutant of $\c{A}$ within $H \cross \c{A}$}
\label{commutant}

\subsection{The basic construction}
\label{basiconstr}

In this subsection we assume that there 
exists an algebra homomorphism
\be
\varphi: H \cross \c{A} \rightarrow \c{A}             \label{Hom}
\ee
acting as the identity on $\c{A}$, namely for any $a\in \c{A}$
\be
\varphi(a)=a .                          \label{ident0}
\ee
(Note that, as a consequence, $\varphi$ is idempotent: 
$\varphi^2=\varphi$). More explicitly,
the fact that $\varphi$ is a homomorphism implies that for any 
$a\in\c{A}$, $g\in H$
\be
\varphi(g)a  =(g_{(1)}\tr a)  \varphi(g_{(2)} ) .                          
\label{fgxrel}
\ee
[Note that (unless the left action of $H$ on $\c{A}$  
is trivial), no $\varphi$ can exist if $\c{A}$ is abelian, 
e.g. for the Euclidean algebra $Uso(N)\cross\b{R}^N$; but,
as we shall recall later, $\varphi$ 
exists for the $q$-deformed Euclidean algebra, and 
therefore the result will apply.]

Let $\c{C}$ be the commutant of $\c{A}$ within $H \cross \c{A}$, 
i.e. the subalgebra 
\be
\c{C}:=\{c\in H \cross \c{A} \: |\:
[c,a]=0  \: \:\: \: \forall a\in\c{A}\}.                \label{commu}
\ee
Clearly $\c{C}$ contains the center $\c{Z}(\c{A})$ of $\c{A}$.
Let $\zeta: H\rightarrow H \cross \c{A}$ be the map defined by
\be
\zeta(g):=\varphi(S g_{(1)}) g_{(2)}.               \label{defzeta}
\ee
Note that if in the definition of $\zeta$ we drop $\varphi$, 
we get instead the counit  $\varepsilon$. Similarly,
if we apply $\varphi$ to $\zeta$ and recall that $\varphi$ is both
a homomorphism and idempotent we also find
\be
\varphi\circ\zeta=\varepsilon.                      \label{billo}
\ee
Now, $\varepsilon(g)$ is a complex number times $\1_{\c{A}}$
and therefore trivially commutes with $\c{A}$. Since
$\varphi$ does not change the commutation relations between $\c{A}$ and 
$ H$, we expect that
also  $\zeta(g)$ commutes with $\c{A}$. This is confirmed by

\begin{theorem} Let $H$ be a Hopf algebra, $\c{A}$ a $H$-module algebra,
$\c{C}$ the commutant (\ref{commu}),
and $\varphi$ an homomorphism of the type (\ref{Hom}), (\ref{ident0}).  Then
(\ref{defzeta}) defines an injective algebra homomorphism
$\zeta: H\rightarrow\c{C}$; moreover $\c{C}=\c{Z}(\c{A})\,\zeta(H)$ 
and $H \cross \c{A}=\c{A}\,\zeta(H)$.
If, in particular  $\c{Z}(\c{A})=\b{C}$, then $\c{C}=\zeta(H)$
and $\zeta: H\leftrightarrow\c{C}$ is an algebra isomorphism.
\label{theorem1}
\end{theorem}

In other words, the subalgebra $H'$ looked for in the introduction
will be obtained by setting $H':=\zeta(H)$.
For these reasons we shall call $\zeta$, as well as the other
maps $\zeta_i,  \zeta_i^{\pm}$ which we shall introduce below,
{\it decoupling maps}.

\bp{}  For any $a\in\c{A}$,
\ba
\zeta(g) a &\stackrel{(\ref{defzeta})}{=}& \varphi(S g_{(1)})g_{(2)} a \nn
&\stackrel{(\ref{gxrel})}{=}&  \varphi(S g_{(1)})(g_{(2)}\tr a)\: g_{(3)} \nn
&\stackrel{(\ref{fgxrel})}{=}& [Sg_{(2)}\tr(g_{(3)}\tr a)] \:
\varphi(S g_{(1)})g_{(4)} \nn
&\stackrel{(\ref{modalg1r})}{=}& [(Sg_{(2)}g_{(3)})\tr a]\:
\varphi(S g_{(1)})g_{(4)}  \nn
&=&[\varepsilon(g_{(2)}) {\bf 1}_{ H} \tr a]\:\varphi(S g_{(1)})g_{(3)} \nn
&=& a \varphi(S g_{(1)})g_{(2)}= a\zeta(g),          \label{comm}
\ea
proving that $\zeta(g)\in\c{C}$.
Using (\ref{comm}) with $a=\varphi(S g'_{(1')})$ we find
\ba
\zeta(gg')&\stackrel{(\ref{defzeta})}{=}&
\varphi[S (g_{(1)}g'_{(1')})]g_{(2)}g'_{(2')}
=\varphi(Sg'_{(1')}) \varphi(S g_{(1)})  g_{(2)}g'_{(2')} \nn
&\stackrel{(\ref{defzeta})}{=}& \varphi(Sg'_{(1')}) \zeta(g) g'_{(2')}
\stackrel{(\ref{comm})}{=}\zeta(g) \varphi(Sg'_{(1')}) g'_{(2')}  \nn
&\stackrel{(\ref{defzeta})}{=}&\zeta(g)\zeta(g'), \label{homomo}
\ea
proving that $\zeta$ is a homomorphism. To prove that $\zeta$ is 
injective note that $\zeta(g)=\zeta(g')$ implies
$$
\varphi(S g_{(1)})\otimes g_{(2)}=\varphi(Sg'_{(1)}) \otimes g'_{(2)},
$$
whence, by applying 
$(m\otimes\id)\circ(\id\otimes\varphi\otimes\id) \circ(\id\otimes\Delta)$
we find $g=g'$ (we have denoted by
$m$ the multiplication map of $\c{A}$, $m(a\otimes b)=ab$).

Now, consider a generic element $c\in H \cross \c{A}$
and decompose it in the form $c=c^{(1)}c^{(2)}$ with 
$c^{(1)}\otimes c^{(2)}\in \c{A}\otimes H $. From (\ref{defzeta})
it immediately follows that 
$$
c=c^{(1)}c^{(2)}=c^{(1)}\varphi\left(c^{(2)}_{(1)}\right)
\zeta\left(c^{(2)}_{(2)}\right),
$$
showing that $H \cross \c{A}=\c{A}\,\zeta(H)$, because
$c^{(1)}\varphi\left(c^{(2)}_{(1)}\right)\in\c{A}$. In particular, assume
$c\in\c{C}$. Then
$$
0=[a,c]=\left[a,c^{(1)}\varphi\left(c^{(2)}_{(1)}\right)\zeta
\left(c^{(2)}_{(2)}\right)\right]\stackrel{(\ref{comm})}{=}
\left[a,c^{(1)}\varphi\left(c^{(2)}_{(1)}\right)\right]
\zeta\left(c^{(2)}_{(2)}\right);
$$
since $\zeta$ is injective, all factors
$\zeta\left(c^{(2)}_{(2)}\right)$ are linearly independent  and therefore
$$
\left[a,c^{(1)}\varphi\left(c^{(2)}_{(1)}\right)\right]=0,
$$
whence
we conclude that $c^{(1)}\varphi\left(c^{(2)}_{(1)}\right)\in\c{Z}(\c{A})$
and $c\in\c{Z}(\c{A})\zeta(H)$. \ep

\begin{corollary} Under the same assumptions of Thm. \ref{theorem1}
the center of the cross-product $H \cross \c{A}$ is given by
\be
\c{Z}(H \cross \c{A})=\c{Z}(\c{A})\,\zeta\left(\c{Z}(H)\right).
\ee
Moreover, if $H_c,\c{A}_c$ are maximal abelian subalgebras of $H$
and $\c{A}$ respectively, then $\c{A}_c\,\zeta(H_c)$
is a maximal abelian subalgebra of $H \cross \c{A}$.
\label{corol1}
\end{corollary}
This is almost all what we need in the determination of 
the Casimirs and of
a complete set of commuting observables of a quantum
system whose algebra of observables is equal to (or contains) 
$H \cross \c{A}$, as in Ref. \cite{Fio95,CerMadSchWes00}.
In addition we just need that these two subalgebras
be closed under the corresponding $*$-structure of 
$H\cross \c{A}$, what will be investigated in section \ref{star}.

\bp{} The proof of the second statement is immediate. As for the
first, if $c\in\c{Z}(H)$ then $[c,H]=0$ and, applying the
homomorphism $\zeta$,
$[\zeta(c),\zeta(H)]=0$; on the other hand $[\zeta(c),\c{A}]=0$
by (\ref{comm}). Hence, $\zeta(c)$ commutes with $\c{A}\zeta(H)$,
i.e. with $H \cross \c{A}$, by theorem \ref{theorem1}.
Similarly if $\tilde c\in\c{Z}(\c{A})$ then $[\tilde c,\c{A}]=0$;
on the other hand $[\tilde c,\zeta(H)]=0$ by (\ref{comm}).
Hence $\tilde c$  commutes with $H \cross \c{A}$. 
Therefore 
$\c{Z}(\c{A})\,\zeta\left(\c{Z}(H)\right)\subset\c{Z}(H\cross\c{A})$.

Viceversa, by theorem \ref{theorem1} any $c\in H\cross\c{A}$ can be
expressed in the form $c=c^{(1)}\zeta\left(c^{(2)}\right)$ with 
$c^{(1)}\otimes c^{(2)}\in \c{A}\otimes H$. If in particular
$c\in\c{Z}(H\cross\c{A})$, then it must be on one hand
$$
0=[c^{(1)}\zeta\left(c^{(2)}\right),\c{A}]=
[c^{(1)},\c{A}]\,\zeta\left(c^{(2)}\right),
$$
implying $c^{(1)}\in\c{Z}(\c{A})$ by the linear independence of all factors 
$\zeta\left(c^{(2)}\right)$; on the other hand it must be
$$
0=[c^{(1)}\zeta\left(c^{(2)}\right),\zeta(H)]=
c^{(1)}[\zeta\left(c^{(2)}\right),\zeta(H)]=
c^{(1)}\zeta\left([c^{(2)},H]\right),
$$
implying $c^{(2)}\in\c{Z}(H)$, by the linear independence of all factors 
$c^{(1)}$ and the injectivity of $\zeta$. Therefore
$\c{Z}(H\cross\c{A})\subset\c{Z}(\c{A})\zeta\left(\c{Z}(H)\right)$.
\ep

Using the results of the theorem we easily show that the restrictions
of the left action of 
$H$ to $\varphi(H)$ and $H$ itself [see (\ref{adjo})] look the same:

\begin{prop} Under the same assumptions of Thm. \ref{theorem1},
on the images of $\varphi$ the left action reads
\be
g\tr \varphi(h)=\varphi(g\tr h);
\ee
equivalently,
\be
g\varphi(h)=\varphi(g_{(1)}\tr h) g_{(2)}.            \label{fggrel}
\ee
\label{prop1}
\end{prop}

\bp{} 
\ba
\varphi(g\tr h) & \stackrel{(\ref{adjo})}{=} &
\varphi(g_{(1)}h S g_{(2)}) =\varphi(g_{(1)})\varphi(h )\varphi(S g_{(2)})\nn
& \stackrel{(\ref{defzeta})}{=} &
\varphi(g_{(1)})\varphi(h )\zeta(g_{(2)})Sg_{(3)}
\stackrel{(\ref{comm})}{=}
\varphi(g_{(1)})\zeta(g_{(2)})\varphi(h)Sg_{(3)}\nn
& \stackrel{(\ref{defzeta})}{=} &g_{(1)}\varphi(h)Sg_{(2)}
\stackrel{(\ref{adjo})}{=}g\tr\varphi(h). \nonumber 
\ea
\ep

Apart from $\zeta_1\equiv\zeta$,
other maps fulfilling the same property (\ref{billo}) are 
\be
\begin{array}{ll}
 & \zeta_2(g):=g_{(2)}\varphi(S^{-1}g_{(1)}) \\
\zeta_3(g):= Sg_{(1)}\varphi(g_{(2)})\qquad
&\zeta_4(g):=\varphi(g_{(2)})S^{-1}g_{(1)}\\
\zeta_5(g):=g_{(1)}\varphi(Sg_{(2)})\qquad
&\zeta_6(g):=\varphi(S^{-1}g_{(2)})g_{(1)} \\
\zeta_7(g):=S^{-1}g_{(2)}\varphi(g_{(1)})\qquad
&\zeta_8(g):=\varphi(g_{(1)})Sg_{(2)} 
\end{array}
\label{defzeta'}
\ee
One could wonder whether they also fulfill the previous
theorems.
Using (\ref{fggrel}) one can easily show that:
\begin{itemize} 
\item 
\be
\zeta_2=\zeta;                                         \label{zz}
\ee 

\item $\zeta_3=\zeta\circ S$, 
$\zeta_4=\zeta_2\circ S^{-1}=\zeta\circ S^{-1}$ so that
$\zeta_3(g),\zeta_4(g)\in\c{C}$, but $\zeta_3,\zeta_4$ are 
antihomomorphisms;

\item $\zeta_5,\zeta_6,\zeta_7,\zeta_8$ do {\it not} map $ H$ into
$\c{C}$.

\end{itemize}

Let us prove for instance the first statement:
\ba
\zeta_2(g)&\stackrel{(\ref{defzeta'})}{=}
&g_{(2)}\varphi(S^{-1}g_{(1)})\stackrel{(\ref{fggrel})}{=}
\varphi(g_{(2)}\tr S^{-1}g_{(1)})g_{(3)} \nn
&\stackrel{(\ref{adjo})}{=}&
\varphi(g_{(2)} S^{-1}g_{(1)}Sg_{(3)}) g_{(4)}
=\varphi(Sg_{(1)}) g_{(2)}
\stackrel{(\ref{defzeta})}{=}\zeta(g). \nonumber
\ea

How does $\zeta(H)$ transform under the action $\tr$ of $H$?
One can easily verify that it is not mapped into itself.
On the contrary, under the opposite action $\tro$ it is:

\begin{prop} Under the same assumptions of Thm. \ref{theorem1},
for any $g,h\in H$
\ba
&&h\tro\zeta(g)=\zeta(h\tro g), \label{prima}\\
&&h\zeta(g)= \zeta(h_{(2)}\tro g)\, h_{(1)}. \label{seconda}
\ea
\label{prop2}
\end{prop}
\bp{}  By (\ref{realiz1}),
(\ref{seconda}) is a direct consequence of (\ref{prima}). To prove 
the latter,
\ba
\zeta(h\tro g)&\stackrel{(\ref{realizop})}{=}&
\zeta(h_{(2)}\,g\,S^{-1}h_{(1)})\stackrel{(\ref{homomo})}{=}
\zeta(h_{(2)})\zeta(g)\zeta(S^{-1}h_{(1)})\nn
&\stackrel{(\ref{defzeta}), (\ref{zz})}{=}&
\zeta_2(h_{(3)})\zeta(g)\varphi(h_{(2)})S^{-1}h_{(1)}\nn
&\stackrel{(\ref{comm})}{=}&
\zeta_2(h_{(3)})\varphi(h_{(2)})\zeta(g)S^{-1}h_{(1)} \nn
&\stackrel{(\ref{defzeta'})}{=}&
h_{(4)}\varphi(S^{-1}h_{(3)})\varphi(h_{2)})\zeta(g)S^{-1}h_{(1)}\nn
&=&h_{(2)}\zeta(g)S^{-1}h_{(1)}\stackrel{(\ref{realizop})}{=}
h\tro\zeta(g). \nonumber
\ea
\ep

\subsection{Construction adapted to Gauss decompositions of $H$}
\label{gaussconstr}
In view of the applications that we shall consider
in section \ref{applications} it is now useful to
consider the case that, instead of a $\varphi$ we just have 
at our disposal two homomorphisms $\varphi^+,\varphi^-$ 
\be
\varphi^{\pm}: H^{\pm} \cross \c{A} \rightarrow \c{A}      \label{Hom+-}
\ee
fulfilling (\ref{ident0}), 
where $H^+,H^-$ denote two Hopf subalgebras of $H$ such that
Gauss  decompositions $H=H^+H^-=H^-H^+$ hold.
(The typical case is when $H=\uqg$ and $H^+,H^-$ denote its
positive and negative Borel subalgebras.)
Then the theorems listed so far will apply separately to 
$H^+ \cross \c{A}$ and $H^- \cross \c{A}$, if we define
corresponding maps $\zeta^{\pm}: H^{\pm} \rightarrow \c{A}$ by
\be
\zeta^{\pm}(g):=\varphi^{\pm}(S g_{(1)}) g_{(2)},       \label{defzeta+-}
\ee
where $g\in H^{\pm}$ respectively. What can we say about
the whole $H\cross \c{A}$? We now prove

\begin{theorem} Let $H$ be a Hopf algebra, $\c{A}$ a $H$-module algebra,
$\c{C}$ the commutant (\ref{commu}),
and $\varphi^{\pm}$ homomorphisms of the type (\ref{Hom+-}), (\ref{ident0}). 
Under the above assumptions formulae
(\ref{defzeta+-}) define injective algebra homomorphisms
$\zeta^{\pm}: H^{\pm}\rightarrow\c{C}$. Moreover,
\be 
\c{C}=\c{Z}(\c{A})\,\zeta^+(H^+)\,\zeta^-(H^-)
=\c{Z}(\c{A})\,\zeta^-(H^-)\,\zeta^+(H^+)
\ee 
and
\be
H \cross \c{A}=\c{A}\,\zeta^+(H^+)\,\zeta^-(H^-)
=\c{A}\,\zeta^-(H^-)\,\zeta^+(H^+).
\ee
In particular, if  $\c{Z}(\c{A})=\b{C}$, 
then $\c{C}=\zeta^+(H^+)\,\zeta^-(H^-)=\zeta^-(H^-)\,\zeta^+(H^+)$.
\end{theorem}
\label{theorem+-}

\bp{} As anticipated,
the fact that $\zeta^{\pm}$ are injective algebra homomorphisms
$\zeta^{\pm}: H^{\pm}\rightarrow\c{C}$ follows from theorem 
\ref{theorem1}. 
Now, by the Gauss  decomposition $H=H^+H^-$
a generic element $c\in H \cross \c{A}$
can be decomposed in the form $c=c^{(1)}c^{(2)}c^{(3)}$ with 
$c^{(1)}\otimes c^{(2)}\otimes c^{(3)}\in \c{A}\otimes H^+ \otimes H^-$;
again, the Sweedler-type notation at the rhs means that a sum of many
terms is understood. From (\ref{defzeta+-})
it immediately follows that 
\ba
c &=&c^{(1)}c^{(2)}c^{(3)} \nn
&=&c^{(1)}\varphi^+\left(c^{(2)}_{(1)}\right)\zeta^+\left(c^{(2)}_{(2)}\right)
\varphi^-\left(c^{(3)}_{(1')}\right)\zeta^-\left(c^{(3)}_{(2')}\right)\nn
&\stackrel{(\ref{comm})}{=}&\zeta^+\left(c^{(2)}_{(2)}\right)
\zeta^-\left(c^{(3)}_{(2')}\right) c^{(1)}\varphi^+\left(c^{(2)}_{(1)}
\right)\varphi^-\left(c^{(3)}_{(1')}\right), \nonumber
\ea
showing that $H \cross \c{A}=\c{A}\zeta^+(H^+)\zeta^-(H^-)$. 
In particular, assume $c\in\c{C}$. Then
\ba
0&=&[a,c]=\left[a,c^{(1)}\varphi^+\left(c^{(2)}_{(1)}
\right)\varphi^-\left(c^{(3)}_{(1')}\right)\zeta^+\left(c^{(2)}_{(2)}\right)
\zeta^-\left(c^{(3)}_{(2')}\right)\right]\nn
&\stackrel{(\ref{comm})}{=}&
\left[a,c^{(1)}\varphi^+\left(c^{(2)}_{(1)}
\right)\varphi^-\left(c^{(3)}_{(1')}\right)\right]
\zeta^+\left(c^{(2)}_{(2)}\right)\zeta^-\left(c^{(3)}_{(2')}\right);\nonumber
\ea
since $\zeta^+,\zeta^-$ are injective, all factors
$\zeta^+\left(c^{(2)}_{(2)}\right)\zeta^-\left(c^{(3)}_{(2')}\right)$ 
are linearly independent  and therefore
$$
\left[a,c^{(1)}\varphi^+\left(c^{(2)}_{(1)}
\right)\varphi^-\left(c^{(3)}_{(1')}\right)\right]=0,
$$
whence
we conclude that $c^{(1)}\varphi^+\left(c^{(2)}_{(1)}
\right)\varphi^-\left(c^{(3)}_{(1')}\right)\in\c{Z}(\c{A})$
and $c$ belongs to
$\c{Z}(\c{A})\,\zeta^+(H^+)\,\zeta^-(H^-)$. 

The proof
of the claims with $\zeta^+(H^+),\zeta^-(H^-)$ in the inverse order
follow in the same way from the Gauss  decomposition $H=H^-H^+$. \ep

As a consequence of this theorem, for any $g^+\in H^+$, $g^-\in H^-$
there exists a sum $c^{(1)}\otimes c^{(2)}\otimes c^{(3)}\in 
\c{Z}(\c{A})\otimes H^-\otimes H^+$
(depending on $g^+,g^-$) such that
\be
\zeta^+(g^+)\zeta^-(g^-)=c^{(1)}\zeta^-(c^{(2)})\zeta^+(c^{(3)}).
                                                    \label{com+-}
\ee
These will be the ``commutation relations'' between elements of
$\zeta^+(H^+)$ and $\zeta^-(H^-)$. Their form will depend on
the specific algebras considered. In section
\ref{applications} we shall determine these commutation
relations for two examples of cross products with $H=\uqg$ using
the Faddeev-Reshetikhin-Takhtadjan generators of $\uqg$.

Of course propositions \ref{prop1}, \ref{prop2} will still apply
to the present situation if both $g,h$ belong to $H^+$
and we replace $\varphi,\zeta$, by $\varphi^+,\zeta^+$
(or the same with $+$ replaced by $-$). What can we say
otherwise?

\begin{prop} Under the same assumptions of Thm. \ref{theorem+-},
if $g\in H^{\pm}$, $h\in H^{\mp}$ 
\ba
&&g\tr \varphi^{\mp}(h)=\varphi^{\pm}(g_{(1)})\varphi^{\mp}(h)\varphi^{\pm}
(Sg_{(2)}),\\
&&g \varphi^{\mp}(h)=\varphi^{\pm}(g_{(1)})\varphi^{\mp}(h)
\varphi^{\pm}(Sg_{(2)})g_{(3)},\\
&&g\tro \zeta^{\mp}(h)=\zeta^{\pm}(g_{(2)})\zeta^{\mp}(h)
\zeta^{\pm}(S^{-1}g_{(1)}),\\
&&g \zeta^{\mp}(h)=\zeta^{\pm}(g_{(3)})\zeta^{\mp}(h)
\zeta^{\pm}(S^{-1}g_{(2)})g_{(1)}.
\ea
\end{prop}
The proof uses theorem \ref{theorem+-} and is similar to the proofs
of Propositions \ref{prop1}, \ref{prop2}. Similarly,
the analog of Corollary \ref{corol1} reads

\begin{corollary} Under the same assumptions of Thm. \ref{theorem+-},
any element $c$ of the center $\c{Z}(H \cross \c{A})$
of the cross-product $H \cross \c{A}$
can be expressed in the form 
\be
c=\zeta^+\left(c^{(1)}\right)\zeta^-\left(c^{(2)}\right)c^{(3)},
\ee
where $c^{(1)}\otimes c^{(2)}\otimes c^{(3)}
\in H^+\otimes H^- \otimes\c{Z}(\c{A})$ and
$c^{(1)}c^{(2)}\otimes c^{(3)}
\in \c{Z}(H) \otimes\c{Z}(\c{A})$; viceversa any
such object $c$ is an element of $\c{Z}(H \cross \c{A})$.
If $H_c\subset H^+ \cap H^-$ and 
$\c{A}_c$ are  maximal abelian subalgebras of $H$
and $\c{A}$ respectively, then $\c{A}_c\,\zeta^+(H_c)$
(as well as $\c{A}_c\,\zeta^-(H_c)$)
is a maximal abelian subalgebra of $H \cross \c{A}$.
\label{corol2}
\end{corollary}

\sect{$*$-structures}
\label{star}

Assume that $H$ is a Hopf $*$-algebra and $\c{A}$ a $H$-module
$*$-algebra, which means that on $H$ and 
$\c{A}$ there exist antilinear involutive antihomomorphisms,
which we both denote by the symbol $*$, such that
\ba
&& [\varepsilon(g)]^*=\varepsilon(g^*)\label{cou*} \\
&& (g^*)_{(1)}\botimes (g^*)_{(2)}=(g_{(1)})^*\botimes (g_{(2)})^*\label{cop*}\\
&& (Sg)^*=S^{-1}g^*  \label{S*}\\
&& (g \tr a_i)^*=(S^{-1} g^*)\tr a_i^* \label{act*}
\ea
[actually (\ref{S*}) is a consequence of  (\ref{cou*}), (\ref{cop*})
and of the uniqueness of the antipode].
Then these two $*$-structures can be glued in a unique one to make 
$H \cross \c{A}$ a $*$-algebra itself. If 
\be
\varphi: H \cross \c{A} \rightarrow \c{A}            
\ee
is a homomorphism acting as the identity on $\c{A}$, then it
is straightforward to check that
$\varphi':= *\circ\varphi\circ *$ also is. As a consequence,
if in a concrete case we know that such a homomorphism is unique, 
then $\varphi'=\varphi$ and we automatically conclude that it is
a $*$-homomorphism,
\be
\varphi(\alpha^*)=[\varphi(\alpha)]^*,\qquad\qquad \alpha\in
H \cross \c{A}.                                 \label{*hom}
\ee
More generally, if the homomorphism $\varphi$ is not
unique, it is natural to look for
one that is a $*$-homomorphism. We shall
give explicit examples of this in section \ref{applications}.

\begin{prop}
If $\varphi: H \cross \c{A} \rightarrow \c{A}$
is a $*$-homomorphism, then also the map
$\zeta: H \rightarrow \c{C}$ is.
\end{prop}

\bp{}
\ba
\zeta(g^*)&\stackrel{(\ref{defzeta})}{=}&\varphi(Sg^*_{(1)})g^*_{(2)}
\stackrel{(\ref{S*}),(\ref{cop*})}{=}
\varphi\left((S^{-1}g_{(1)})^*\right)g^*_{(2)}  \nn
&\stackrel{(\ref{*hom})}{=}&\left[\varphi\left(S^{-1}g_{(1)}\right)\right]^*
g^*_{(2)}=\left[g_{(2)}\varphi\left(S^{-1}g_{(1)}\right)\right]^* \nn
&\stackrel{(\ref{defzeta'})}{=}&[\zeta_2(g)]^*\stackrel{(\ref{zz})}{=}
[\zeta(g)]^*.         \nonumber
\ea
\ep

Alternatively, it may happen that no $*$-homomorphism
$\varphi$ exists, but
there exist homomorphisms $\varphi^{\pm}$ 
of the type (\ref{Hom+-}). If $H^{\pm}$ are Hopf 
$*$-subalgebras (i.e. are closed under $*$), then also
$\varphi'{}^{\pm}:= *\circ\varphi^{\pm}\circ *$  are  
homomorphisms of the type (\ref{Hom+-}), and as before
we can look for $\varphi^{\pm}$ that are $*$-homomorphisms,
\be
\varphi^{\pm}(\alpha^*)=[\varphi^{\pm}(\alpha)]^*,\qquad\qquad \alpha\in
H^{\pm} \cross \c{A}.                                \label{*hom+-}
\ee
On the contrary, if $H^{\pm}$  are mapped into each other by $*$,
then $*\circ\varphi^{\pm}\circ *$ is a homomorphism of the
type $\varphi^{\mp}$, and we can look for one such that
\be
\varphi^{\pm}(\alpha^*)=[\varphi^{\mp}(\alpha)]^*,\qquad\qquad \alpha\in
H^{\mp} \cross \c{A}.                                \label{*hom+-'}
\ee

\begin{prop}
If $\varphi^{\pm}$ are $*$-homomorphism, then also the map
$\zeta^{\pm}: H^{\pm} \rightarrow \c{C}$ are. If $\varphi^{\pm}$
fulfill (\ref{*hom+-'}), then $\zeta^{\pm}$ fulfill
\be
\zeta^{\pm}(g^*)=[\zeta^{\mp}(g)]^*,\qquad\qquad g\in
H^{\mp}.                              
\ee
\end{prop}

\bp{} The first statement amounts to the preceding
proposition applied to $H^{\pm} \cross \c{A}$. As for the second,
the proof is just a small variation:
\ba
\zeta^{\pm}(g^*)&\stackrel{(\ref{defzeta})}{=}&\varphi^{\pm}
(Sg^*_{(1)})g^*_{(2)}\stackrel{(\ref{S*}),(\ref{cop*})}{=}
\varphi^{\pm}\left((S^{-1}g_{(1)})^*\right)g^*_{(2)}  \nn
&\stackrel{(\ref{*hom+-'})}{=}&\left[\varphi^{\mp}
\left(S^{-1}g_{(1)}\right)\right]^*
g^*_{(2)}=\left[g_{(2)}\varphi^{\mp}\left(S^{-1}g_{(1)}\right)\right]^* \nn
&\stackrel{(\ref{defzeta'})}{=}&[\zeta^{\mp}_2(g)]^*\stackrel{(\ref{zz})}{=}
[\zeta^{\mp}(g)]^*.         \nonumber
\ea
\ep

\sect{Formulae for right-module algebras}
\label{rightsec}

In this section we give the analogs for right $H$-module algebras
of the main results found so far for left $H$-module algebras.
By definition the right action $\tl: \c{A} \times H\to \c{A}$ fulfills
\ba
&&a\tl(gg') = (a\tl g) \tl   g'                     \label{modalg1}\\
&&(aa')\tl g = (a\tl g_{(1)})\, (a'\tl g_{(2)}).       \label{modalg2}
\ea
The algebra $\c{A} \cocross H$ is defined
as follows. 
As a vector space it is $H\otimes\c{A}$, whereas
the product is defined through formulae obtained from
(\ref{trivial1}-\ref{trivial3}) by flipping the tensor
factors, together with
$$
(\1_H\otimes a)(g\otimes\1_{\c{A}} )=g_{(1)} \otimes (a\tl g_{(2)})
$$
for any $a\in\c{A}$, $g\in H$. As before, we shall denote 
$g\otimes a$ by $ga$ and omit either unit $\1_{\c{A}},\1_H$ 
whenever multiplied
by non-unit elements; consequently only the last condition
takes a non-trivial form and becomes
\be
a g=g_{(1)}\, (a\tl g_{(2)}).                        \label{crossprod}
\ee
$\c{A}\cocross H$ itself is a $H$-module algebra under the right
action  $\tl$ of $ H$ if we extend the latter on the elements of $H$
as the adjoint action,
\be
h\tl g=Sg_{(1)} h g_{(2)},          \qquad\qquad g,h\in  H; \label{adjor}
\ee
Let $\tilde\c{C}$ be the commutant of $\c{A}$ within $\c{A}\cocross H$.
Clearly it contains the center $\c{Z}(\c{A})$ of $\c{A}$.
We shall need homomorphisms
\be
\tilde\varphi: \c{A}\cocross H \rightarrow \c{A}   \label{Homr}
\ee
acting as the identity on $\c{A}$, namely for any $a\in \c{A}$
\be
\tilde\varphi(a)=a.                                 \label{ident0r}
\ee

\begin{theorem} Let $H$ be a Hopf algebra, $\c{A}$ a right
$H$-module algebra, $\tilde\c{C}$ the commutant  
of $\c{A}$ within $\c{A}\cocross H$,
and $\tilde\varphi$ an homomorphism of the type (\ref{Homr}), 
(\ref{ident0r}). Then the map $\zeta_5$ defined by
(\ref{defzeta'}) (with $\varphi$ replaced by $\tilde\varphi$)
is an injective algebra homomorphism
$\zeta_5: H\rightarrow\tilde\c{C}$; moreover 
$\tilde\c{C}=\zeta_5(H)\,\c{Z}(\c{A})$ and
$\c{A}\cocross H=\zeta_5(H)\,\c{A}$.
If, in particular  $\c{Z}(\c{A})=\b{C}$, then $\tilde\c{C}=\zeta(H)$
and $\zeta_5: H\leftrightarrow\tilde\c{C}$ is an algebra isomorphism.
\end{theorem}
\label{theorem1r}

Consider the maps $\zeta_i$ defined by
(\ref{defzeta'}), but with $\varphi$ replaced by $\tilde\varphi$. 
One can easily show that:
\begin{itemize} 
\item 
\be
\zeta_6=\zeta_5;                                         \label{zzr}
\ee 

\item $\zeta_8=\zeta_6\circ S=\zeta_5\circ S$, 
$\zeta_7=\zeta_5\circ S^{-1}$, so that
$\zeta_7(g),\zeta_8(g)\in\c{C}$, but $\zeta_7,\zeta_8$ are 
antihomomorphisms;

\item $\zeta_1,\zeta_2,\zeta_3,\zeta_4$ do {\it not} map $H$ into
$\tilde\c{C}$.

\end{itemize}
\begin{prop} If $g,h\in H$
\ba
&&\tilde\varphi(h)\tl g=\tilde\varphi(h\tl g), \\
&&\tilde\varphi(h)g=g_{(1)}\tilde\varphi(h\tl g_{(2)}).    \label{fgrel}
\ea
\end{prop}

If, instead of a $\tilde\varphi$, we just have 
at our disposal two homomorphisms $\tilde\varphi^+,\tilde\varphi^-$ 
\be
\tilde\varphi^{\pm}: \c{A}\cocross H^{\pm}  \rightarrow \c{A},   \label{Hom+-r}
\ee
where $H^+,H^-$ denote two Hopf subalgebras of $H$ such that
Gauss  decompositions $H=H^+H^-=H^-H^+$ hold,
then the theorems listed so far will apply separately to 
$\c{A}\cocross H^+$ and $\c{A}\cocross H^-$, if we define
corresponding maps $\zeta_5^{\pm}: H^{\pm} \rightarrow \c{A}$ by
\be
\zeta_5^{\pm}(g):=g_{(1)}\tilde\varphi^{\pm}(Sg_{(2)}),   \label{defzeta+-r}
\ee
where $g\in H^{\pm}$ respectively. More precisely:

\begin{theorem} Under the above assumptions formulae
(\ref{defzeta+-r}) define injective algebra homomorphisms
$\zeta_5^{\pm}: H^{\pm}\rightarrow\tilde\c{C}$. Moreover,
\be 
\tilde\c{C}=\zeta_5^+(H^+)\,\zeta_5^-(H^-)\,\c{Z}(\c{A})
=\zeta_5^-(H^-)\,\zeta_5^+(H^+)\,  \c{Z}(\c{A})       \label{casim} 
\ee 
and
\be
\c{A} \cocross H=\zeta_5^+(H^+)\,\zeta_5^-(H^-)\,\c{A}
=\zeta_5^-(H^-)\,\zeta_5^+(H^+)\,\c{A}.         \label{decom}
\ee
In particular, if  $\c{Z}(\c{A})=\b{C}$, 
then $\tilde\c{C}=\zeta_5^+(H^+)\zeta_5^-(H^-)=\zeta_5^-(H^-)\zeta_5^+(H^+)$.
\label{theorem+-r}
\end{theorem}

\begin{prop}
If $g\in H^{\pm}$, $h\in H^{\mp}$ 
\ba
&&\tilde\varphi^{\mp}(h)\tl g=\tilde\varphi^{\pm}(Sg_{(1)})\tilde
\varphi^{\mp}(h)\tilde\varphi^{\pm}(g_{(2)}),\\
&& \tilde\varphi^{\mp}(h) g=g_{(1)}\tilde\varphi^{\pm}(Sg_{(2)})\tilde
\varphi^{\mp}(h)\tilde\varphi^{\pm}(g_{(3)}).
\ea
\end{prop}

If $H$ is a Hopf $*$-algebra and $\c{A}$ a $H$-module
$*$-algebra, then 
the two $*$-structures of $H$ and $\c{A}$
can be glued into a unique one to make 
$\c{A}\cocross H$ a $*$-algebra itself.

\begin{prop}
If $\tilde\varphi: \c{A}\cocross H\rightarrow \c{A}$
is a $*$-homomorphism, then also the map
$\zeta_5: H \rightarrow \tilde\c{C}$ is.
\label{zeta5*}
\end{prop}

\begin{prop}
If the maps $\tilde\varphi^{\pm}$ defined
in (\ref{defzeta+-r}) are $*$-homomorphism, then also the maps
$\zeta_5^{\pm}: H^{\pm} \rightarrow \c{C}$ are. If $\tilde\varphi^{\pm}$
fulfill
\be
\tilde\varphi^{\pm}(\alpha^*)=[\tilde\varphi^{\mp}(\alpha)]^*,
\qquad\qquad \alpha\in
\c{A}\cocross H^{\mp},                                \label{*hom+-'r}
\ee
then $\zeta_5^{\pm}$ fulfill
\be
\zeta_5^{\pm}(g^*)=[\zeta_5^{\mp}(g)]^*,\qquad\qquad g\in
H^{\mp}.                              
\ee
\label{zeta5+-*}
\end{prop}

\sect{Applications}
\label{applications}

We now consider a couple of applications where
$H$ is the quantum group \uqg \cite{Dri86}, with
$\g=sl(N)$ or $\g=so(N)$.
As a set of generators of \uqg
it is convenient to introduce the FRT generators \cite{FadResTak89}
$\c{L}^+{}_j^i,\c{L}^-{}_j^i$ ($i,j$ take
$N$ different values), together with
the square roots of the
diagonal elements $\c{L}^+{}_j^j,\c{L}^-{}_j^j$.
In the appendix we recall the relations they fulfill.
The FRT generators are related to the so-called
universal $R$-matrix  $\R$ by
\be
\c{L}^+{}_l^a:=\R^{(1)}\rho_l^a(\R^{(2)})\qquad\qquad
\c{L}^-{}_l^a:=\rho_l^a(\R^{-1}{}^{(1)})\R^{-1}{}^{(2)}, \label{frt}
\ee
where we have denoted by $\rho$ the fundamental $N$-dimensional
representation of $U_qsl(N)$ or $U_qso(N)$.
Since in our conventions $\R\in H^+\botimes H^-$ 
($H^+,H^-$ denote the positive, negative Borel subalgebras) we see that 
$\c{L}^+{}_l^a\in H^+$ and $\c{L}^-{}_l^a\in H^-$. 

For historical reasons we introduce algebras $\c{A}$ 
as right- (rather than left-) $\uqg$-module algebras.

\subsection{The Euclidean quantum group $\b{R}_q^N\cocross U_qso(N)$}
\label{appliEuc}

As algebra $\c{A}$ we shall consider a slight extension of
the quantum Euclidean space $\b{R}_q^N$
~\cite{FadResTak89} (the $U_qso(N)$-covariant quantum space),
i.e.  of the unital associative algebra generated by $p^i$ fulfilling the
relations 
\be 
\c{P}_a{}^{ij}_{hk}p^hp^k=0, \label{xxrel} 
\ee 
where $\c{P}_a$
denotes the $q$-deformed antisymmetric projector appearing in the decomposition
of the braid matrix $\hat R$ of $U_qso(N)$ [given in formula (\ref{defRsoN})];
the latter is related to $\R$ by 
$\hat R^{ij}_{hk}=\rho^j_h(\c{L}^+{}_k^i)=(\rho^j_h\otimes\rho^i_k)(\R)$.
The multiplet $(p^i)$ carries the fundamental $N$-dim
(or vector) representation $\rho$
of $U_qso(N)$: for any $g\in U_qso(N)$
\be 
p^i\tl g=\rho^i_j(g)p^j.  \label{fund1} 
\ee
This implies
\ba
&&p^i\c{L}^{\pm}{}^a_b=\c{L}^{\pm}{}^a_c p^j\hat 
R^{\pm 1}{}^{ci}_{jb},                     \label{gio}\\
&&S\c{L}^{\pm}{}^a_b p^i=\hat 
R^{\pm 1}{}^{ai}_{jk}p^jS\c{L}^{\pm}{}^k_b. \label{giu}
\ea
To define $\tilde\varphi$ or 
$\tilde\varphi^{\pm}$ one \cite{CerFioMad00} slightly enlarges 
$\b{R}_q^N$ as follows. One 
introduces some new generators $\sqrt{P_a}$, with
$1\le a\le \frac N 2$, together with their  
inverses $(\sqrt{P_a})^{-1}$, requiring that 
\be
P_a^2=\sum\limits_{h=-a}^a p^hp_h=\sum\limits_{h,k=-a}^a g_{hk}p^hp^k.
\ee
In the previous equation $g_{hk}$ denotes the `metric matrix' of 
$SO_q(N)$:
\be
g_{ij}=g^{ij}=q^{-\rho_i} \delta_{i,-j}.          \label{defgij}
\ee
It is a $SO_q(N)$-isotropic tensor and is a deformation of the 
ordinary Euclidean metric.
Here and in the sequel $n:=\left[\frac N 2\right]$ is the rank of 
$so(N)$, the indices take the values
$i=-n,\ldots,-1,0,1,\ldots n$ for $N$ odd,
and $i=-n,\ldots,-1, 1,\ldots n$ for $N$ even.
We have also introduced the notation
$(\rho_i)=(n-\frac{1}{2},\ldots,\frac{1}{2},0,-\frac{1}{2},
\ldots,\frac{1}{2}-n)$
for $N$ odd, $(n-1,\ldots,0,0,\ldots,1-n)$ for $N$ even. 
In the sequel we shall call $P_n^2$ also $P^2$.
Moreover for odd $N$ we add also $\sqrt{p^0}$ and its inverse
as new generators. The commutation relations 
involving these new generators can be fixed consistently, and turn out to
be simply $q$-commutation relations.
$P$ plays the role of `deformed Euclidean distance' 
of the generic `point of coordinates' $(p^i)$ of $\b{R}_q^N$ 
from the `origin';
$P_a$ is the `projection' of $P$ on the `subspace' $p^i=0$, $|i|>a$. 
The center of $\b{R}_q^N$ is generated by $\sqrt{P}$ and, only in the case
of even $N$, by $\sqrt{\frac{p^1}{p^{-1}}}$ and its inverse
$\sqrt{\frac{p^{-1}}{p^1}}$.
In the case of even $N$ one needs to include also
the FRT generator  $\c{L}^-{}^1_1=\c{L}^+{}^{-1}_{-1}$ 
and its inverse $\c{L}^+{}^1_1=\c{L}^-{}^{-1}_{-1}$
[which are generators of $U_qso(N)$ belonging to the natural Cartan 
subalgebra] among the generators of $\c{A}$. From (\ref{gio}) one
derives that  they satisfy the commutation relations
\be
\c{L}^-{}^1_1 p^{\pm 1}=q^{\pm 1} p^{\pm 1} \c{L}^-{}^1_1, \quad\quad
\c{L}^-{}^1_1 p^{\pm i}= p^{\pm i} \c{L}^-{}^1_1\:\: \hbox{ for } i>1
\label{xkapparel}
\ee
with the generators of $\c{A}$, and the standard FRT relations with the rest
of $U_qso(N)$. As a consequence, 
$\sqrt{\frac{p^{\pm 1}}{p^{\mp1}}}$ are eliminated from the center of
$\c{A}$ (in fact $\c{L}^{\pm}{}^1_1$
do not $q$-commute with $\sqrt{\frac{p^{\pm 1}}{p^{\mp1}}}$).

One can easily show that the extension of the action of 
$U_qso(N)$
to $\sqrt{P_a},(\sqrt{P_a})^{-1}$ is uniquely determined by the constraints
the latter fulfill; it is a bit complicated and therefore will 
be omitted, since we will not need its explicit expression. 
We keep the action of $H$ on $\c{L}^-{}^1_1$ as the standard 
(right) adjoint action (\ref{adjor}). Note that the maps $\tilde\varphi^{\pm}$ 
have no analog
in the ``undeformed'' case ($q=1$),
because $\c{A}_1\equiv\b{R}^N$ is abelian, whereas 
$H\equiv U_qso(N)$ is not. 

The homomorphisms \cite{CerFioMad00} 
$\tilde\varphi^{\pm}:\c{A} \cocross  U_q^{\pm}so(N) \rightarrow \c{A}$
take the simplest and most compact expression on the FRT generators of
$U_q^{\pm}so(N)$. Let us introduce the short-hand notation
$[A,B]_x=AB-xBA$.     
The images of $\tilde\varphi^-$ on the negative FRT generators
read
\be
\tilde\varphi^-(\c{L}^-{}^i_j)=g^{ih}[\mu_h,p^k]_qg_{kj},      \label{imagel-}
\ee
where
\be
\begin{array}{ll}
\mu_0=\gamma_0 (p^0)^{-1}
&\quad\mbox{for $N$ odd,} \\[6pt]
\mu_{\pm 1}=\gamma_{\pm 1} (p^{\pm 1})^{-1} \c{L}^{\pm}{}^1_1
&\quad\mbox{for $N$ even,} \\[6pt]
\mu_a=\gamma_a P_{|a|}^{-1}P_{|a|-1}^{-1} p^{-a}
&\quad\mbox{otherwise,} 
\end{array}                                             \label{defmu}
\ee
and $\gamma_a \in \b{C}$ are normalization constants fulfilling
the conditions
\be
\begin{array}{ll}
\gamma_0 = -q^{-\frac{1}{2}} h^{-1} &\quad\mbox{for $N$ odd,} \\[6pt]
\gamma_{\pm 1}=-k^{-1} &\quad\mbox{for $N$ even,}\\[6pt]
\gamma_1 \gamma_{-1}=-q^{-1} h^{-2} &\quad\mbox{for $N$ odd,} \\[6pt]
\gamma_a \gamma_{-a} =
-q^{-1} k^{-2} \omega_a \omega_{a-1} &\quad\mbox{for $a>1$}. \nonumber
\end{array}                                               \label{gamma}
\ee
Here $h:=q^{\frac{1}{2}}-q^{-\frac{1}{2}}$,
$k:=q\!-\!q^{-1}$, $\omega_a:=(q^{\rho_a}+q^{-\rho_a})$.
On the other hand, the images of $\tilde\varphi^+$ on the positive 
FRT generators read
\be
\tilde\varphi^+(\c{L}^+{}^i_j)=g^{ih}[\bar\mu_h,p^k]_{q^{-1}}g_{kj},
                                                     \label{imagel+}
\ee
where
\be
\begin{array}{ll}
\bar\mu_0=\bar\gamma_0 (p^0)^{-1}
&\quad\mbox{for $N$ odd,} \\[6pt]
\bar\mu_{\pm 1} = 
\bar\gamma_{\pm 1} (p^{\pm 1})^{-1} \c{L}^{\mp}{}^1_1
&\quad\mbox{for 
$N$ even,} \\[6pt]
\bar\mu_a = 
\bar\gamma_a P_{|a|}^{-1}P_{|a|-1}^{-1} p^{-a}
&\quad\mbox{otherwise,} 
\end{array}                                         \label{defbarmu}
\ee
and $\bar\gamma_a \in \b{C}$ normalization constants fulfilling
the conditions
\be
\begin{array}{ll}
\bar \gamma_0 = q^{\frac{1}{2}} h^{-1}&\quad\mbox{for $N$ odd,} \\[6pt]
\bar \gamma_{\pm 1}=k^{-1}&\quad\mbox{for $N$ even,}\\[6pt]
\bar \gamma_1 \bar \gamma_{-1} = -q h^{-2}&\quad\mbox{for $N$ odd,} \\[6pt]
\bar \gamma_a \bar \gamma_{-a} = -q k^{-2} \omega_a \omega_{a-1}
&\quad\mbox{for $a>1$}.         \nonumber
\end{array}                                             \label{bargamma}
\ee

From definition (\ref{defzeta'}), 
using (\ref{coprodL}), (\ref{imagel-}), (\ref{imagel+}), 
(\ref{antip}), we find
\be
\begin{array}{l}
\zeta_5^-(\c{L}^-{}^i_j)=\c{L}^-{}^i_h\tilde\varphi^-(S\c{L}^-{}^h_j)=
\c{L}^-{}^i_h[\mu_h,p^i]_q, \\
\zeta_5^+(\c{L}^+{}^i_j)=\c{L}^+{}^i_h\tilde\varphi^+(S\c{L}^+{}^h_j)=
\c{L}^+{}^i_h[\bar\mu_h,p^i]_{q^{-1}}. 
\end{array}          \label{z+-}
\ee

\medskip
The case of even $N$ is slightly outside the scheme developed in the
preceding sections. As anticipated, the generators $\c{L}^{\pm}{}^1_1$
and their inverses $\c{L}^{\pm}{}^{-1}_{-1}$ cannot be realized
as (rational) `functions' of the $p$'s and have to be introduced
in the codomain of $\tilde\varphi^{\pm}$ as new generators. 
In fact definitions (\ref{imagel-}), (\ref{imagel+}) are exactly designed
to lead to
\be
\tilde\varphi^{\pm}(\c{L}^{\pm}{}^1_1)=\c{L}^{\pm}{}^1_1
\qquad\qquad
\tilde\varphi^{\pm}(\c{L}^{\pm}{}^{-1}_{-1})=\c{L}^{\pm}{}^{-1}_{-1}.
\label{ciccio}
\ee
As a consequence, when $s=\pm 1$
\ba
&&\zeta_5^+(\c{L}^+{}^s_s)=\c{L}^+{}^s_k\tilde\varphi^+(S\c{L}^+{}^k_s)
\stackrel{(\ref{sfilza1}),(\ref{ciccio})}{=} \c{L}^+{}^s_sS\c{L}^+{}^s_s=1\\
&&\zeta_5^-(\c{L}^-{}^s_s)=\c{L}^-{}^s_k\tilde\varphi^-(S\c{L}^-{}^k_s)
\stackrel{(\ref{sfilza2}),(\ref{ciccio})}{=} \c{L}^-{}^s_sS\c{L}^-{}^s_s=1.
\ea
Moreover, it is easy to check that
\be
\begin{array}{l}
\c{L}^-{}^1_1\zeta_5^{\pm}(\c{L}^{\pm}{}^i_j)=
\zeta_5^{\pm}(\c{L}^{\pm}{}^i_j)\c{L}^-{}^1_1 q^{\eta_i-\eta_j}\\
\c{L}^+{}^1_1\zeta_5^{\pm}(\c{L}^{\pm}{}^i_j)=
\zeta_5^{\pm}(\c{L}^{\pm}{}^i_j)\c{L}^+{}^1_1 q^{\eta_j-\eta_i},
\end{array}
\label{tete'}
\ee
where we have introduced the shorthand notation
\be
\eta_i:=\delta_i^1-\delta_i^{-1}.                   \label{short}
\ee
Therefore, $\zeta_5^{\pm}[U^{\pm}_qso({2n})]$ do not commute
with the whole $\c{A}$ but only with $\b{R}_q^{2n}$,
and the decomposition (\ref{casim}) 
will have to be modified into
\be 
\tilde\c{C}=\zeta_5^+(H^+)\zeta_5^-(H^-)\c{Z}(\b{R}_q^{2n})
=\zeta_5^-(H^-)\zeta_5^+(H^+)\c{Z}(\b{R}_q^{2n})   \label{casim'} 
\ee 
with $\tilde\c{C}$ defined as the commutant of $\b{R}_q^{2n}$ within 
$\b{R}_q^{2n}\cocross U_qso({2n})$. As said, $\c{Z}(\b{R}_q^{2n})$ contains
also $\sqrt{\frac{p^{\pm 1}}{p^{\mp1}}}$.
Since $\mu_{\pm 1}$, $\bar\mu_{\mp 1}$
and therefore $\tilde\varphi^-(\c{L}^-{}^i_{\pm 1})$,
$\tilde\varphi^+(\c{L}^+{}^i_{\mp 1})$ are of the form
``$\c{L}^{\mp}{}^1_1~\times$  an expression depending only on the 
$p$'s'', they $q$-commute rather than commute with
$\zeta^{\pm}(\c{L}^{\pm}{}^i_j)$. 

\medskip

For odd and even $N$
the commutation relations among the $\zeta_5^-(\c{L}^-{}^i_j)$
or among the $\zeta_5^+(\c{L}^+{}^i_j)$ are immediately obtained
from (\ref{sfilza1}-\ref{L-L-rel}) and (\ref{LLg})
applying $\zeta_5^-,\zeta_5^+$ and 
using the fact that they are homomorphisms. 

To derive the commutation relations between
the $\zeta_5^-(\c{L}^-{}^i_j)$
and the $\zeta_5^+(\c{L}^+{}^h_k)$
[the analog of (\ref{com+-}) in explicit form] 
we need the commutation relations
between the $\tilde\varphi^-(\c{L}^-{}^i_j)$
and the $\tilde\varphi^+(\c{L}^+{}^i_j)$. 
To proceed we have to distinguish the case of odd and even $N$.
In the formulation of the following lemma
and proposition we switch off Einstein summation convention.
In the appendix we prove

\begin{lemma} If $N$ is odd, then
\be
\tilde\varphi^-(S\c{L}^-{}^i_k)\tilde\varphi^+(S\c{L}^+{}^h_j)=
\frac{\gamma_k}{\bar\gamma_k}\sum\limits_{l,m,r,s}
\hat R^{-1}{}^{ih}_{lm}\tilde\varphi^+(S\c{L}^+{}^l_r)
\tilde\varphi^-(S\c{L}^-{}^m_s)\hat R^{rs}_{kj}
\frac{\bar\gamma_s}{\gamma_s}.                   \label{phi+-}
\ee
If $N$ is even, then if $k\notin \{-1,1\}$
\ba
&&\tilde\varphi^-(S\c{L}^-{}^i_k)\tilde\varphi^+(S\c{L}^+{}^h_j) =
\frac{\gamma_k}{\bar\gamma_k}\sum\limits_{l,m}\hat R^{-1}{}^{ih}_{lm}
\left[\sum\limits_{r,s\atop s\neq \pm 1}\tilde\varphi^+(S\c{L}^+{}^l_r)
\tilde\varphi^-(S\c{L}^-{}^m_s)\hat R^{rs}_{kj}
\frac{\bar\gamma_s}{\gamma_s}\right. \nn
&&\qquad\qquad\qquad\left.+\sum\limits_{r,s\atop s= \pm 1}
\tilde\varphi^+(S\c{L}^+{}^l_r)
\tilde\varphi^-(S\c{L}^-{}^{~m}_{-s})\frac{p^{-s}}{p^s}\hat R^{rs}_{kj}
\frac{\bar\gamma_s}{\gamma_s}q^2\right],       \label{phi+-'}
\ea
and if $k\in \{-1,1\}$
\ba
&&\tilde\varphi^-(S\c{L}^-{}^i_k)\tilde\varphi^+(S\c{L}^+{}^h_j) =
-\!\sum\limits_{l,m}\hat R^{-1}{}^{ih}_{lm}\!
\left[\!\sum\limits_{r,s\atop s\neq \pm 1}\tilde\varphi^+(S\c{L}^+{}^l_r)
\tilde\varphi^-(S\c{L}^-{}^m_s)\hat R^{~~rs}_{-k,j}
\frac{\bar\gamma_s}{\gamma_s}\right. \nn
&&\qquad\left.+\sum\limits_{r,s\atop s= \pm 1}
\tilde\varphi^+(S\c{L}^+{}^l_r)
\tilde\varphi^-(S\c{L}^-{}^{~m}_{-s})\frac{p^{-s}}{p^s}\hat R^{~~rs}_{-k,j}
\frac{\bar\gamma_s}{\gamma_s}q^2\right]\frac{p^{-k}}{p^k}q^{-2+2\eta_j\eta_k}.
\label{phi+-"}
\ea
\label{lala}
\end{lemma}

Note that for odd
$N$ if the ratio $\frac{\gamma_a}{\bar\gamma_a}$ is independent
of $a$ (therefore it must be equal to
$\frac{\gamma_0}{\bar\gamma_0}=-q^{-1}$), 
the $\gamma,\bar\gamma$'s disappear,
and by comparison with (\ref{L+L-rel}) we find \cite{CerFioMad00}
that setting
$\tilde\varphi(\c{L}^+{}^i_j)=\tilde\varphi^+(\c{L}^+{}^i_j)$,
$\tilde\varphi(\c{L}^-{}^i_j)=\tilde\varphi^-(\c{L}^-{}^i_j)$
defines a homomorphism $\tilde\varphi:\c{A}\cocross U_qso(N)\to\c{A}$.
For $|q|=1$ $\tilde\varphi$ turns out to be a $*$-homomorphism
w.r.t. the corresponding $*$-structures (see below).

\begin{prop} If $N$ is odd
\be
\zeta_5^+(\c{L}^+{}^h_j)\zeta_5^-(\c{L}^-{}^i_k)=
\sum\limits_{c,d,r,s}
\hat R^{-1}{}^{ih}_{cd}\zeta_5^-(\c{L}^-{}^d_s)
\zeta_5^+(\c{L}^+{}^c_r)\hat R^{rs}_{kj}
\frac{\gamma_k}{\bar\gamma_k}\frac{\bar\gamma_s}{\gamma_s}.
                                             \label{z+-rel}
\ee
If $N$ is even, then
\ba
&&\zeta_5^+(\c{L}^+{}^h_j)\zeta_5^-(\c{L}^-{}^i_k)=
\frac{\gamma_k}{\bar\gamma_k}q^{\eta_j(\eta_i-\eta_k)}
\sum\limits_{l,m,r}\hat R^{-1}{}^{ih}_{lm}
\left[\sum\limits_{s\neq \pm 1}
\zeta_5^-(\c{L}^-{}^m_s)\zeta_5^+(\c{L}^+{}^l_r)\right.\nn
&&\left.+\sum\limits_{s= \pm 1}
\zeta_5^-(\c{L}^-{}^m_{-s})\zeta_5^+(\c{L}^+{}^l_r)
\frac{p^{-s}}{p^s}q^{2+\eta_s(\eta_r-\eta_l)}\right]
\hat R^{rs}_{kj}
\frac{\bar\gamma_s}{\gamma_s}               \label{z+-rel'}
\ea
if $k\notin\{1,-1\}$, and, if $k\in\{1,-1\}$,
\ba
&&\zeta_5^+(\c{L}^+{}^h_j)\zeta_5^-(\c{L}^-{}^i_k)=
-q^{-2+\eta_j(\eta_i+\eta_k)}\sum\limits_{l,m,r}\hat R^{-1}{}^{ih}_{lm}
\left[\sum\limits_{s\neq \pm 1}
\zeta_5^-(\c{L}^-{}^m_s)\zeta_5^+(\c{L}^+{}^l_r)\right.\nn
&&\left.+\sum\limits_{s= \pm 1}
\zeta_5^-(\c{L}^-{}^m_{-s})\zeta_5^+(\c{L}^+{}^l_r)\frac{p^{-s}}{p^s}
q^{2+\eta_s(\eta_r-\eta_l)}\right]
\zeta_5^+(\c{L}^+{}^l_r)\hat R^{~~rs}_{-k,j}
\frac{\bar\gamma_s}{\gamma_s}\frac{p^{-k}}{p^k}.       \label{z+-rel"}
\ea
\label{pippo}
\end{prop}
So for even $N$ the coefficients in the commutations relations depend
explicitly on the central (in $\b{R}_q^N$) elements 
$\frac{p^{\pm 1}}{p^{\mp1}}$.

Let us analyze now the properties of 
$\zeta_5^{\pm}$ under $*$-structures.
When $|q|=1$ the $*$-structure of $\b{R}_q^N$
is given by $(p^i)^*=p^i$. 
It turns out \cite{FioSteWes00}
that $\tilde\varphi^{\pm}$ are $*$-homomorphisms if, in addition, 
\be
\gamma_a^*=-\gamma_a\left\{
\begin{array}{lll}
1 &\mbox{ if }\:\:\:a<-1,\qquad &\mbox{ or } a=-1 \mbox{ and } N\mbox{ odd;}\\
q^{-2}&\mbox{ if }\:\:\: a>1,\qquad &\mbox{ or } a=1  \mbox{ and } N\mbox{ odd.} 
\end{array}\right. 
\label{REAL}
\ee
Under these assumptions, in view also of (\ref{|q|=1}), we can apply 
proposition
\ref{zeta5+-*} and we conclude that $\zeta_5^{\pm}$ are $*$-homomorphisms.
When $q\in\b{R}^+$ the real structure of $\b{R}_q^N$ is given by
$(p^i)^*=p^jg_{ji}$.
It turns out \cite{FioSteWes00} that $\tilde\varphi^{\pm}$ 
fulfill $[\varphi^{\pm}(g)]^*=\varphi^{\mp}(g^*)$, or more
explicitly
\be
\left[\varphi^-(\c{L}^-{}^i_j)\right]^* =
\varphi^+\left[(\c{L}^-{}^i_j)^*\right],
\ee
if, in addition, 
\be
\gamma_a^*=-\bar\gamma_{-a}\left\{
\begin{array}{lll}
1 &\mbox{ if }\:\:\:a<-1,\qquad &\mbox{ or } a=-1 
\mbox{ and } N\mbox{ odd;}\\
q^{-2}&\mbox{ if }\:\:\: a>1,\qquad &\mbox{ or } a=1  
\mbox{ and } N\mbox{ odd.} 
\end{array}\right. 
\ee
Under these assumptions, in view also of (\ref{qReal}),
(\ref{|q|=1}), we can apply proposition
\ref{zeta5+-*} and we conclude that $\zeta_5^{\pm}$ 
fulfill
\be
\zeta_5^{\pm}(g^*)=[\zeta_5^{\mp}(g)]^*,\qquad\qquad g\in
H^{\mp}.                              
\ee

\subsection{The cross product of \uqg's with  \uqg-covariant 
Heisenberg algebras}
\label{heisenberg}

As algebra $\c{A}$ we shall consider a slight extension of the
$\uqg$-covariant deformed Heisenberg algebras 
$\c{D}_{\epsilon,\g}$, $\g=sl(N),so(N)$. Such algebras 
have been introduced in Ref. 
\cite{PusWor89,WesZum90,CarSchWat91}. They are
unital associative algebras generated by $x^i,\partial_j$
fulfilling the relations
\begin{equation}
\c{P}_a{}^{ij}_{hk}x^hx^k=0, \qquad
\c{P}_a{}^{ij}_{hk}\partial_j\partial_i=0,
\qquad\partial_i x^j=\delta^i_j+(q\gamma\hat R)^{\epsilon}
{}_{ih}^{jk}x^h\partial_k,
\label{gringo}
\end{equation}
where $\gamma=q^{\frac 1N},1$ respectively for $\g=sl(N),so(\!N\!)$,
and the exponent $\epsilon$ can take either value $\epsilon=1,-1$.
$\hat R$ denotes the braid matrix of $\uqg$ [given in formulae
(\ref{defRslN}), (\ref{defRsoN})], and the matrix
$\c{P}_a$ is the deformed
antisymmetric projector appearing in the decompositions 
(\ref{projectorR}), (\ref{projectorR'}) of the latter. 
The coordinates $x^i$, as the $p^i$ of
subsection \ref{appliEuc}, transform according to the fundamental
$N$-dimensional representation $\rho$ of $\uqg$, 
whereas the `partial derivatives' transform according the
contragradient representation,
\begin{equation}
x^i\tl g=\rho^i_j(g)x^j,   \qquad\qquad
\partial_i\tl g=
\partial_h\rho^h_i(S^{-1}g).                  \label{fund}
\end{equation}
In our conventions
the indices will take the values  $i=1,...,N$ if $\g=sl(N)$,
whereas if $\g=so(N)$ they will take the same values considered
in subsection \ref{appliEuc}. In fact the quantum Euclidean 
space can be considered as a subalgebra of the 
$U_qso(N)$-covariant Heisenberg algebra, either
by the identification $p^i\equiv x^i$ or by the one
$p^i\equiv g^{ij}\partial_j$.

Algebra homomorphisms
$\tilde\varphi:\c{A}\cocross H\rightarrow \c{A}$, for 
$H=\uqg$  and $\c{A}$ equal to
(a suitable completion of) $\c{D}_{\epsilon,\g}$
have been constructed in Ref. \cite{Fiocmp95,ChuZum95}.
This is the $q$-analog of the well-known fact that the elements of  
$\g$ can be realized as ``vector fields''
(first order differential operators) on the corresponding
$\g$-covariant (undeformed) space, e.g.
$\tilde\varphi(E^i_j)=x^i\partial_j-\frac 1N \delta^i_j$
in the $\g=sl(N)$ case.
This means that our decoupling map $\zeta_5$ exists
also in the undeformed case; we don't know whether this
result in Lie group representation theory has ever been
formulated before.

The explicit expression of $\tilde\varphi(\c{L}^-{}^i_j)$
in terms of $x^i,\partial_j$ for $U_qsl(2),U_qso(3)$ 
has been given in Ref. \cite{FioSteWes00}. For different values
of $N$ it can be found from the results of Ref. 
\cite{Fiocmp95,ChuZum95} by passing from the generators
adopted there to the FRT generators.
For example, for $\g=sl(2)$ and $\epsilon=1$
one finds 
\ba
&&\tilde\varphi(\c{L}^+{}^1_1)=\tilde\varphi(\c{L}^-{}^2_2)=
[\tilde\varphi(\c{L}^-{}^1_1)]^{-1}=[\tilde\varphi(\c{L}^+{}^2_2)]^{-1}=
\alpha\Lambda^{\frac 12}\left[1\!+\!(q^2\!-\!1)
x^2\partial_2\right]^{\frac 12}\nn
&&\tilde\varphi(\c{L}^+{}^1_2)=-\alpha kq^{-1}\Lambda^{\frac 12}\left[1+(q^2-1)
x^2\partial_2\right]^{-\frac 12}x^1\partial_2 \\
&&\tilde\varphi(\c{L}^-{}^2_1)=\alpha kq^3\Lambda^{\frac 12}\left[1+(q^2-1)
x^2\partial_2\right]^{-\frac 12}x^2\partial_1,\nonumber
\ea
where $\alpha$ is fixed by (\ref{L+L-rel}) to be
$\alpha=\pm 1,\pm i$ and we have set
\be
\Lambda^{-2}:=1+(q^2-1)x^i\partial_i.        \label{defLambda-2'}
\ee
Whereas for $\g=so(3)$ and $\epsilon=1$
one finds on the positive Borel subalgebra
\be
\begin{array}{l}
\tilde\varphi(\c{L}^+{}^-_-)=-\alpha\Lambda\left[1+(q-1)x^0\partial_0+(q^2-1)
x^+\partial_+\right]\\
\tilde\varphi(\c{L}^+{}^-_0)=\alpha k\Lambda(x^-\partial_0-\sqrt{q}x^0\partial_+)\\
\tilde\varphi(\c{L}^+{}^-_+)=\frac 1{1+q^{-1}}\tilde\varphi(\c{L}^+{}^-_0)
\tilde\varphi(\c{L}^+{}^0_+)\\
\tilde\varphi(\c{L}^+{}^0_0)=1 \\
\tilde\varphi(\c{L}^+{}^0_+)=-q^{-\frac 12}[\tilde\varphi(\c{L}^+{}^-_-)]^{-1}
\tilde\varphi(\c{L}^+{}^-_0)\\
\tilde\varphi(\c{L}^+{}^+_+)= [\tilde\varphi(\c{L}^+{}^-_-)]^{-1}
\end{array} \label{part1}
\ee
and on the negative Borel subalgebra
\be
\begin{array}{l}
\tilde\varphi(\c{L}^-{}^-_-)=-\left(\alpha\Lambda\left[1+(q-1)x^0
\partial_0+(q^2-1)x^+\partial_+\right]\right)^{-1}\\
\tilde\varphi(\c{L}^-{}^0_-)=-\alpha q^2 k \tilde\varphi(\c{L}^-{}^-_-)
\Lambda(x^0\partial_--\sqrt{q}x^+\partial_0)\\
\tilde\varphi(\c{L}^-{}^+_-)=\frac1{1+q}\tilde\varphi(\c{L}^-{}^+_0)
\tilde\varphi(\c{L}^-{}^0_-)\\
\tilde\varphi(\c{L}^-{}^0_0)=1 \\
\tilde\varphi(\c{L}^-{}^+_0)=-\alpha q^{\frac 32}k\Lambda(x^0\partial_--
\sqrt{q}x^+\partial_0)\\
\tilde\varphi(\c{L}^-{}^+_+)= [\tilde\varphi(\c{L}^-{}^-_-)]^{-1}. 
\end{array}
\label{part2}
\ee
Here we have set
\be
\Lambda^{-2}:=[1+(q^2-1)x^i\partial_i+\frac{(q^2-1)^2}{\omega_1^2}
(g_{ij}x^ix^j)(g^{hk}\partial_k\partial_h)],        \label{defLambda-2}
\ee
where
$$
\omega_a:=(q^{\rho_a}+q^{-\rho_a}), 
$$
and replaced
for simplicity the values $-1,0,1$ of the indices by the ones $-,0,+$. 
In either case the $\tilde\varphi$-images of $\c{L}^+{}^i_j$ and
$\c{L}^-{}^j_i$ for $i>j$ vanish, because the latter do.

We see that strictly speaking $\tilde\varphi$ takes values in
some appropriate completion of $\c{D}_{\epsilon,\g}$, containing
at least the square root and inverse square root of the polynomial 
$\Lambda^{-2}$ respectively
 defined in (\ref{defLambda-2'}), (\ref{defLambda-2}),
as well as the square root of $[1+(q^2-1)
x^2\partial_2]$ and its inverse,
when $\g=sl(2)$, and the inverses
(\ref{part1})$_6$, (\ref{part2})$_6$, when $\g=so(3)$. 
Apart from this minimal completion, another possible one
is the so-called $h$-adic, namely the ring
of formal power series in $h=\log q$ with coefficients in
$\c{D}_{\epsilon,\g}$. Other completions, e.g. in operator norms,
can be considered according to the needs.
One can easily show that the extension of the action of $H$ 
to any such completion is uniquely determined (we omit to write 
down its explicit expression, since we don't  need it).

$\c{A}\cocross H$ is a $*$-algebra and
the map $\tilde\varphi$ is a $*$-homomorphism both for
$q$ real and $|q|=1$. The $*$-structure of $\c{A}$ is
\be
(x^i)^*=x^i, \qquad (\partial_i)^*=-\partial_i\cases{q^{\pm 2(N-i+1)}
\:\:\mbox{ if }H=U_qsl(N)\cr
q^{\pm N+\rho_i}\:\:\mbox{ if }H=U_qso(N)}
\ee
if $|q|=1$, and 
\be
(x^h)^*=x^kg_{kh}, \qquad (\partial_i)^*=-\frac{\Lambda^{\pm 2}}{
q^{\pm N}+q^{\pm 2}}\left[(g^{jh}\partial_h\partial_j),\, x^i\right]
\ee
if $H=U_qso(N)$ and $q\in\b{R}^+$. The upper or lower sign respectively refer
to the choices
$\epsilon=1,-1$ in (\ref{gringo})$_3$, and 
$\Lambda^{\pm 2}$ are respectively defined by
\be
\Lambda^{\pm 2}:=\left[1+(q^{\pm 2}-1)x^i\partial_i+ \frac{(q^{\pm 2}
-1)^2}{\omega_n^2} r^2 (g^{ji}\partial_i\partial_j)\right]^{-1}.
\ee
In either case, in view also of (\ref{qReal}), we can apply 
proposition
\ref{zeta5*} and we conclude that  $\zeta_5$ is a
$*$-homomorphism.

\app{Appendix}

\subsubsection*{Basic properties of \uqg}

For both $H=U_qsl(N),U_qso(N)$ the FRT generators fulfill the 
relations
\ba
&&\c{L}^+{}^i_j=0,           \hspace{1.5cm}\mbox{if $i>j$},\label{sfilza1}\\
&&\c{L}^-{}^i_j=0,           \hspace{1.5cm}\mbox{if $i<j$},\label{sfilza2}\\
&&\c{L}^{-}{}^i_i\c{L}^{+}{}^i_i=\c{L}^{+}{}^i_i\c{L}^{-}{}^i_i=1,
\hspace{1cm}\forall i  \label{sfilza3}\\
&&\prod_i\c{L}^+{}^i_i=1,\qquad\prod_i\c{L}^-{}^i_i=1,     \label{sfilza4}\\
&&\hat R^{ab}_{cd}\,\c{L}^+{}^d_f\c{L}^+{}^c_e=
\c{L}^+{}^b_c\c{L}^+{}^a_d\,\hat R^{dc}_{ef},      \label{L+L+rel}\\
&&\hat R^{ab}_{cd}\,\c{L}^-{}^d_f\c{L}^-{}^c_e=
\c{L}^-{}^b_c\c{L}^-{}^a_d\,\hat R^{dc}_{ef}, \label{L-L-rel}\\
&&\hat R^{ab}_{cd}\,\c{L}^+{}^d_f\c{L}^-{}^c_e=
\c{L}^-{}^b_c\c{L}^+{}^a_d\,\hat R^{dc}_{ef};\label{L+L-rel}
\ea
in addition, when $H=U_qso(N)$ they also fulfill
\be
\c{L}^{\pm}{}^i_j\c{L}^{\pm}{}^h_k g^{kj}=g^{hi}\qquad 
\c{L}^{\pm}{}_i^j\c{L}^{\pm}{}_h^k g_{kj}=g_{hi}.  \label{LLg}
\ee
Here $\hat R^{ab}_{cd}$ denotes the braid matrix of 
$U_qsl(N),U_qso(N)$, and $g^{hi}$ the metric matrix of 
$U_qso(N)$, which was defined in (\ref{defgij}).
The square roots of the diagonal elements $\c{L}^-{}^i_i$,
as well as the the square roots of the
diagonal elements $\c{L}^+{}^i_i$, generate the same Cartan
subalgebra of $H$, which coincides with $H^+\cap H^-$.
The braid matrix $\hat R$ of $U_qsl(N)$ is given by
\be
\hat R = q^{-\frac 1N}\left[q \sum_i e^i_i \otimes e^i_i +
\sum_{\scriptstyle i \neq j} e^j_i \otimes e^i_j
+k \sum_{i<j} e^i_i \otimes e^j_j \right]            \label{defRslN}
\ee
where all indices $i,j,a,...=1,2,...,N$,
and $e^i_j$ is the $N \times N$ matrix with all elements
equal to zero except for a $1$ in the $i$th column and $j$th row.
When $H=U_qso(N)$ it is convenient to adopt the convention that
all indices $i,j,a,...$ take the values
$i=-n,\ldots,-1,0,1,\ldots n$ for $N$ odd,
and $i=-n,\ldots,-1, 1,\ldots n$ for $N$ even, where
$n:=\left[\frac N 2\right]$ is the rank of $so(N)$. 
Then the corresponding braid matrix reads
\ba
\hat R&=&q \sum_{i \neq 0} e^i_i \otimes e^i_i +
\sum_{\stackrel{\scriptstyle i \neq j,-j} 
{\mbox{ or } i=j=0}} e^j_i \otimes e^i_j+ q^{-1} 
\sum_{i \neq 0} e^{-i}_i
\otimes e^i_{-i} 
\label{defRsoN} \\
&&+k (\sum_{i<j} e^i_i \otimes e^j_j- 
\sum_{i<j} q^{-\rho_i+\rho_j} 
e^{-j}_i \otimes e^j_{-i}). \nonumber
\ea

The braid matrix of $sl(N)$ admits the orthogonal projector
decomposition
\be
 q^{\frac 1N}\hat R = q\c{P}_S - q^{-1}\c{P}_a, \qquad\qquad \g=sl(N);
\label{projectorR}  
\ee
$\c{P}_a,\c{P}_S$ are the $U_qsl(N)$-covariant
deformed antisymmetric and symmetric projectors. 
The braid matrix of $so(N)$ admits the orthogonal projector
decomposition
\be
\hat R = q\c{P}_s - q^{-1}\c{P}_a + q^{1-N}\c{P}_t\qquad\qquad \g=so(N);      
\label{projectorR'}
\ee
$\c{P}_a,\c{P}_t,\c{P}_s$ are the corresponding
$q$-deformed antisymmetric, trace,
trace-free symmetric projectors. 

By iterated use of equations (\ref{L+L+rel}-\ref{L+L-rel})
one immediately shows that they hold 
also if $\hat R$ is replaced by any polynomial function
$f(\hat R)$ of $\hat R$, in
particular by $f(\hat R)=\hat R^{-1},\c{P}_a$.

Finally, the coproduct of the FRT generators is given by
\be
\Delta(\c{L}^+{}^i_j)=\c{L}^+{}^i_h\otimes\c{L}^+{}^h_j, \qquad\qquad
\Delta(\c{L}^-{}^i_j)=\c{L}^-{}^i_h\otimes\c{L}^-{}^h_j.
\label{coprodL}
\ee
When $H=U_qso(N)$ the antipode is given by
\be
S\c{L}^{\mp}{}^j_i=g_{ih}\c{L}^{\mp}{}^h_k g^{kj}.        \label{antip}
\ee

The non-compact real sections of \uqg require $|q|=1$ and are
characterized by the $*$-structure
\be
(\c{L}^{\pm}{}^i_j)^*=U^{-1}{}^i_r\,\c{L}^{\pm}{}^r_s\, U^s_j. \label{|q|=1}
\ee
Here one can take
\be
U^i_j:=\cases{g^{ih}g_{jh} \mbox{ if }\g=so(N),\cr     
               q^{-i}\delta^i_j \mbox{ if }\g=sl(N)}.      \label{defU}
\ee

The compact real section of \uqg requires $q\in\b{R}^+$ if $\g=so(N)$, 
$q\in\b{R}$ if $\g=sl(N)$ and is characterized by the $*$-structure
\be
(\c{L}^{\pm}{}^i_j)^*=S\c{L}^{\mp}{}^j_i.                 \label{qReal}
\ee
For $\g=so(N)$ this amounts to
\be
(\c{L}^{\pm}{}^i_j)^*=g_{ih}\c{L}^{\mp}{}^h_k g^{kj}.        \label{qreal0}
\ee

\subsubsection*{Proof of Lemma \ref{lala}}

We start by recalling
two relations proved in Lemma 2 of Ref. \cite{CerFioMad00}
\be
\begin{array}{l}
\mu_a\tilde\varphi^-(S\c{L}^-{}^i_b) =
\hat{R}^{-1}{}^{cd}_{ab}\tilde\varphi^-(S\c{L}^-{}^i_c) \mu_d,\\
\bar\mu_a \tilde\varphi^+(S\c{L}^+{}^i_b)=
\hat{R}^{cd}_{ab} \tilde\varphi^+(S\c{L}^+{}^i_c)\bar\mu_d.  
\end{array}
\label{Rmux}
\ee
Applying $\tilde\varphi^{\pm}$ to (\ref{giu}) we find
\be
\tilde\varphi^{\pm}(S\c{L}^{\pm}{}^a_b) p^i=
\hat R^{\pm 1}{}^{ai}_{jk}p^j\tilde\varphi^{\pm}(S\c{L}^{\pm}{}^k_b).
\label{lulu}
\ee
We finally note that
\be
\bar\mu_{\pm 1}=-\mu_{\mp 1}\frac{p^{\mp 1}}{p^{\pm 1}}q^2.
\label{mubarmu}
\ee
The claim is a direct consequence of relations
(\ref{Rmux}),  (\ref{lulu}), (\ref{mubarmu}).

\subsubsection*{Proof of Proposition \ref{pippo}}

For odd $N$
\ba
\mbox{lhs}(\ref{z+-rel})
&\stackrel{(\ref{z+-})}{=} &\sum\limits_a
\c{L}^+{}^h_a\tilde\varphi^+(S\c{L}^+{}^a_j)\zeta_5^-(\c{L}^-{}^i_k) \nn
&\stackrel{\mbox{\scriptsize Thm \ref{theorem+-r}}}{=} &\sum\limits_a
\c{L}^+{}^h_a\zeta_5^-(\c{L}^-{}^i_k)\tilde\varphi^+(S\c{L}^+{}^a_j) \nn
&\stackrel{(\ref{z+-})}{=} &\sum\limits_{a,b}
\c{L}^+{}^h_a\c{L}^-{}^i_b\tilde\varphi^-(S\c{L}^-{}^b_k)
\tilde\varphi^+(S\c{L}^+{}^a_j) \nn
&\stackrel{(\ref{phi+-})}{=} &\sum\limits_{a,b}
\sum\limits_{l,m,r,s}\c{L}^+{}^h_a\c{L}^-{}^i_b
\hat R^{-1}{}^{ba}_{lm}\tilde\varphi^+(S\c{L}^+{}^l_r)
\tilde\varphi^-(S\c{L}^-{}^m_s)\hat R^{rs}_{kj}
\frac{\gamma_k}{\bar\gamma_k}\frac{\bar\gamma_s}{\gamma_s} \nn
&\stackrel{(\ref{L+L-rel})}{=} &\sum\limits_{c,d}
\sum\limits_{l,m,r,s}\hat R^{-1}{}^{ih}_{cd}
\c{L}^-{}^d_m\c{L}^+{}^c_l\tilde\varphi^+(S\c{L}^+{}^l_r)
\tilde\varphi^-(S\c{L}^-{}^m_s)\hat R^{rs}_{kj}
\frac{\gamma_k}{\bar\gamma_k}\frac{\bar\gamma_s}{\gamma_s} \nn
&\stackrel{(\ref{z+-})}{=} &\sum\limits_{c,d}
\sum\limits_{m,r,s}\hat R^{-1}{}^{ih}_{cd}
\c{L}^-{}^d_m \zeta_5^+(\c{L}^+{}^c_r)
\tilde\varphi^-(S\c{L}^-{}^m_s)\hat R^{rs}_{kj}
\frac{\gamma_k}{\bar\gamma_k}\frac{\bar\gamma_s}{\gamma_s} \nn
&\stackrel{\mbox{\scriptsize Thm \ref{theorem+-r}}}{=} &\sum\limits_{c,d}
\sum\limits_{m,r,s}\hat R^{-1}{}^{ih}_{cd}
\c{L}^-{}^d_m\tilde\varphi^-(S\c{L}^-{}^m_s) \zeta_5^+(\c{L}^+{}^c_r)
\hat R^{rs}_{kj}
\frac{\gamma_k}{\bar\gamma_k}\frac{\bar\gamma_s}{\gamma_s} \nn
&\stackrel{(\ref{z+-})}{=} &
\mbox{rhs}(\ref{z+-rel}).\nonumber
\ea
For even $N$
\ba
&&\zeta_5^+(\c{L}^+{}^h_j)\zeta_5^-(\c{L}^-{}^i_k)
\stackrel{(\ref{z+-})}{=} \sum\limits_b
\c{L}^+{}^h_b\tilde\varphi^+(S\c{L}^+{}^b_j)\zeta_5^-(\c{L}^-{}^i_k) \nn
&&\stackrel{\mbox{\scriptsize Thm \ref{theorem+-r}},
(\ref{tete'})}{=} 
\sum\limits_b
\c{L}^+{}^h_b\zeta_5^-(\c{L}^-{}^i_k)\tilde\varphi^+(S\c{L}^+{}^b_j) 
q^{\eta_j(\eta_i-\eta_k)}\nn
&&\stackrel{(\ref{z+-})}{=} \sum\limits_{a,b}
\c{L}^+{}^h_b\c{L}^-{}^i_a\tilde\varphi^-(S\c{L}^-{}^a_k)
\varphi^+(S\c{L}^+{}^b_j)q^{\eta_j(\eta_i-\eta_k)}\nn
&&\stackrel{(\ref{L+L-rel})}{=} \sum\limits_{a,b}
\sum\limits_{l,m,c,d}\hat R^{-1}{}^{ih}_{lm}
\c{L}^-{}^m_d\c{L}^+{}^l_c\hat R^{cd}_{ab}
\tilde\varphi^-(S\c{L}^-{}^a_k)
\varphi^+(S\c{L}^+{}^b_j)q^{\eta_j(\eta_i-\eta_k)}.     \nonumber
\ea
If $k\notin\{1,-1\}$ then
\ba
&&\zeta_5^+(\c{L}^+{}^h_j)\zeta_5^-(\c{L}^-{}^i_k)
\stackrel{(\ref{phi+-'})}{=}  \frac{\gamma_k}{\bar\gamma_k}
\sum\limits_{l,m,c,d}\hat R^{-1}{}^{ih}_{lm}
\c{L}^-{}^m_d\c{L}^+{}^l_c
\left[\sum\limits_{r,s\atop s\neq \pm 1}\tilde\varphi^+(S\c{L}^+{}^c_r)
\right. \nn
&&\:\left.\times\tilde\varphi^-(S\c{L}^-{}^d_s)\hat R^{rs}_{kj}
\frac{\bar\gamma_s}{\gamma_s}+\sum\limits_{r,s\atop s= \pm 1}
\tilde\varphi^+(S\c{L}^+{}^c_r)
\tilde\varphi^-(S\c{L}^-{}^d_{-s})\frac{p^{-s}}{p^s}\hat R^{rs}_{kj}
\frac{\bar\gamma_s}{\gamma_s}q^2\right]q^{\eta_j(\eta_i-\eta_k)}\nn
&&\stackrel{(\ref{z+-})}{=} \frac{\gamma_k}{\bar\gamma_k}
\sum\limits_{l,m,d,r}\hat R^{-1}{}^{ih}_{lm}
\c{L}^-{}^m_d\zeta_5^+(\c{L}^+{}^l_r)
\left[\sum\limits_{s\neq \pm 1}
\tilde\varphi^-(S\c{L}^-{}^d_s)\hat R^{rs}_{kj}
\frac{\bar\gamma_s}{\gamma_s}\right. \nn
&&\qquad\left.+\sum\limits_{s= \pm 1}
\tilde\varphi^-(S\c{L}^-{}^d_{-s})\frac{p^{-s}}{p^s}\hat R^{rs}_{kj}
\frac{\bar\gamma_s}{\gamma_s}q^2\right]q^{\eta_j(\eta_i-\eta_k)}\nn
&&\stackrel{\mbox{\scriptsize Thm \ref{theorem+-r}},
(\ref{tete'})}{=} \frac{\gamma_k}{\bar\gamma_k}
\sum\limits_{l,m,d,r}\hat R^{-1}{}^{ih}_{lm}\c{L}^-{}^m_d
\left[\sum\limits_{s\neq \pm 1}
\tilde\varphi^-(S\c{L}^-{}^d_s)\hat R^{rs}_{kj}
\frac{\bar\gamma_s}{\gamma_s}\right. \nn
&&\left. +\sum\limits_{s= \pm 1}
\tilde\varphi^-(S\c{L}^-{}^d_{-s})
\frac{p^{-s}}{p^s}\hat R^{rs}_{kj}
\frac{\bar\gamma_s}{\gamma_s}q^{2+\eta_s(\eta_r-\eta_l)}\right]
\zeta_5^+(\c{L}^+{}^l_r)q^{\eta_j(\eta_i-\eta_k)}\nn
&&\stackrel{(\ref{z+-})}{=} \frac{\gamma_k}{\bar\gamma_k}
\sum\limits_{l,m,r}\hat R^{-1}{}^{ih}_{lm}
\left[\sum\limits_{s\neq \pm 1}
\zeta_5^-(\c{L}^-{}^m_s)\hat R^{rs}_{kj}
\frac{\bar\gamma_s}{\gamma_s}+\right. \nn
&&\left.\sum\limits_{s= \pm 1}
\zeta_5^-(\c{L}^-{}^m_{-s})\frac{p^{-s}}{p^s}\hat R^{rs}_{kj}
\frac{\bar\gamma_s}{\gamma_s}q^{2+\eta_s(\eta_r-\eta_l)}\right]
\zeta_5^+(\c{L}^+{}^l_r)q^{\eta_j(\eta_i-\eta_k)}\nn
&&\stackrel{\mbox{\scriptsize Thm \ref{theorem+-r}}}{=}
\mbox{rhs}(\ref{z+-rel'}).\nonumber
\ea
If $k\in\{1,-1\}$ then
\ba
&&\zeta_5^+(\c{L}^+{}^h_j)\zeta_5^-(\c{L}^-{}^i_k)
\stackrel{(\ref{phi+-'})}{=}  -\!
\sum\limits_{l,m,c \atop d,r}\hat R^{-1}{}^{ih}_{lm}
\c{L}^-{}^m_d\c{L}^+{}^l_c\tilde\varphi^+(S\c{L}^+{}^c_r)
\!\left[\sum\limits_{s\neq \pm 1}
\tilde\varphi^-(S\c{L}^-{}^d_s)\right. \nn
&&\left.\times\frac{\bar\gamma_s}{\gamma_s}
\hat R^{~~rs}_{-k,j}+\sum\limits_{s\atop s= \pm 1}
\tilde\varphi^-(S\c{L}^-{}^d_{-s})\frac{p^{-s}}{p^s}\hat R^{rs}_{-kj}
\frac{\bar\gamma_s}{\gamma_s}q^2\right]q^{\eta_j(\eta_i-\eta_k)}
\frac{p^{-k}}{p^k}q^{-2+2\eta_j\eta_k}\nn
&&\stackrel{(\ref{z+-})}{=} -
\sum\limits_{l,m,d,r}\hat R^{-1}{}^{ih}_{lm}
\c{L}^-{}^m_d\zeta_5^+(\c{L}^+{}^l_r)
\left[\sum\limits_{s\neq \pm 1}
\tilde\varphi^-(S\c{L}^-{}^d_s)\hat R^{~~rs}_{-k,j}
\frac{\bar\gamma_s}{\gamma_s}\right. \nn
&&\qquad\left.+\sum\limits_{s= \pm 1}
\tilde\varphi^-(S\c{L}^-{}^d_{-s})\frac{p^{-s}}{p^s}\hat R^{~~rs}_{-k,j}
\frac{\bar\gamma_s}{\gamma_s}q^2\right]
q^{-2+\eta_j(\eta_i+\eta_k)}\frac{p^{-k}}{p^k}\nn
&&\stackrel{\mbox{\scriptsize Thm \ref{theorem+-r}}, 
(\ref{tete'})}{=} -
\sum\limits_{l,m,r}\hat R^{-1}{}^{ih}_{lm}\c{L}^-{}^m_d
\left[\sum\limits_{s\neq \pm 1}
\tilde\varphi^-(S\c{L}^-{}^d_s)\hat R^{~~rs}_{-k,j}
\frac{\bar\gamma_s}{\gamma_s}\right. \nn
&&\left. +\sum\limits_{s= \pm 1}
\tilde\varphi^-(S\c{L}^-{}^d_{-s})
\frac{p^{-s}}{p^s}\hat R^{~~rs}_{-k,j}
\frac{\bar\gamma_s}{\gamma_s}q^{2+\eta_s(\eta_r-\eta_l)}\right]
\zeta_5^+(\c{L}^+{}^l_r)q^{-2+\eta_j(\eta_i+\eta_k)}\frac{p^{-k}}{p^k}\nn
&&\stackrel{(\ref{z+-})}{=} -
\sum\limits_{l,m,r}\hat R^{-1}{}^{ih}_{lm}
\left[\sum\limits_{s\neq \pm 1}
\zeta_5^-(\c{L}^-{}^m_s)\hat R^{~~rs}_{-k,j}
\frac{\bar\gamma_s}{\gamma_s}+\sum\limits_{s= \pm 1}
\zeta_5^-(\c{L}^-{}^m_{-s})\right. \nn
&&\qquad\left.\times\frac{p^{-s}}{p^s}\hat R^{~~rs}_{-k,j}
\frac{\bar\gamma_s}{\gamma_s}q^{2+\eta_s(\eta_r-\eta_l)}\right]
\zeta_5^+(\c{L}^+{}^l_r)q^{-2+\eta_j(\eta_i+\eta_k)}\frac{p^{-k}}{p^k}\nn
&&\stackrel{\mbox{\scriptsize Thm \ref{theorem+-r}}}{=}
\mbox{rhs}(\ref{z+-rel"}). \nonumber
\ea


\begin{thebibliography}{99}

\bibitem{CarSchWat91} 
U.~Carow-Watamura, M.~Schlieker, S.~Watamura, 
``$SO_q(N)$ covariant differential calculus on quantum space 
and quantum deformation of {S}chroedinger equation'', 
{\em Z.~Physik~C - Particles and Fields} {\bf 49} (1991), 439.


\bibitem{CerFioMad00}
B. L. Cerchiai, G. Fiore, J. Madore,
`` Geometrical Tools for Quantum Euclidean Spaces'',
{\em Commun. Math. Phys.} {\bf 217} (2001), 521-554.

\bibitem{CerMadSchWes00}
B.L. Cerchiai, J. Madore, S. Schraml, J. Wess,
``Structure of the Three-dimensional Quantum Euclidean Space'',
{\em Eur. Phys. J.} {\bf C 16} (2000), 169.

\bibitem{ChuZum95}
C.-S. Chu, B. Zumino,
``Realization of vector fields for quantum groups as pseudodifferential
operators on quantum spaces'', Proc.  XX
Int. Conf. on Group Theory Methods in Physics, Toyonaka (Japan), 1995,
and q-alg/9502005.

\bibitem{Dri86}
V.~Drinfeld, ``Quantum groups,'' in {\em {I.C.M.} Proceedings, {B}erkeley},
  p.~798.
\newblock 1986.

\bibitem{FadResTak89}
L.D.~Faddeev, N.Y.~Reshetikhin, L.~Takhtadjan, 
``Quantization of {L}ie groups and Lie algebras'', 
{\em Alge. i Analy.} {\bf 1} (1989), 178, translated from the
Russian in {\em Leningrad Math. J.} {\bf 1} (1990), 193.

\bibitem{Fiocmp95} G. Fiore,
``Realization of $U_q(so(N))$ within the Differential Algebra
on ${\bf R}_q^N$'',

\bibitem{Fio95} G.\ Fiore, 
``The Euclidean Hopf algebra $U_q(e^N)$ and its fundamental Hilbert space
representations'',
{\em J. Math. Phys.} {\bf 36} (1995), 4363-4405;
``The $q$-Euclidean Algebra $U_q(e^N)$ and
the Corresponding $q$-Euclidean Lattice'', 
{\em Int. J. Mod. Phys.} {\bf A11} (1996), 863-886.

\bibitem{FioMad99}
G. Fiore e J. Madore
``The geometry of the quantum Euclidean space''
{\em J. Geom. Phys.} {\bf 33} (2000), 257-287. 

\bibitem{FioSteWes00}
G. Fiore, H. Steinacker and J. Wess
``Unbraiding the braided tensor product'',  math/0007174.

\bibitem{Ogi92}
O.\ Ogievetsky
``Differential operators on quantum spaces for $GL_q(n)$ and $SO_q(n)$''
{\em Lett. Math. Phys.} {\bf 24} (1992), 245.

\bibitem{fuzzyq} H. Grosse, J. Madore, H. Steinacker,
  ``Field Theory on the $q$--deformed Fuzzy Sphere'', 
{\em J. Geom. Phys.} {\bf 38} (2001), 308-342.

\bibitem{Hay90}
T.\ Hayashi, ``q-Analogs of Clifford and Weyl Algebras: Spinor and
Oscillator Realizations of Quantum Enveloping algebras''
Commun. Math. Phys. {\bf 127} (1990), 129.

\bibitem{OgiSchWesZum92}
O.\ Ogievetsky, W.B.\ Schmidke, J.\ Wess, B.\ Zumino,
``$q$-deformed Poincar\'e algebra'' {\em Commun. Math. Phys.}
{\bf 150} (1992), 495.

\bibitem{PusWor89} W. \ Pusz, S. \ L. \ Woronowicz, 
``Twisted Second Quantization'',
{\em Rep. Math. Phys.} {\bf 27} (1989), 231.

\bibitem{WesZum90}
J.~Wess, B.~Zumino, 
``Covariant differential calculus on the quantum hyperplane'', 
{\em Nucl.\ Phys.\ (Proc.\ Suppl.)} {\bf 18B} (1990) 302.

\end{thebibliography}
\end{document}